\documentclass[a4paper,11pt,reqno]{amsart}

\usepackage[utf8]{inputenc}

\usepackage{amsmath}
\usepackage{amsthm}
\usepackage{amssymb}
\usepackage{mathtools}

\usepackage{tikz-cd}

\usepackage[all,cmtip]{xy}
\usepackage{hyperref}
\usepackage[capitalize]{cleveref}
\usepackage{bm}
\usepackage{mathdots}
\usepackage{array}
\usepackage{dsfont}
\usepackage{xcolor}
\usepackage{pdfpages}
\usepackage[
hscale=0.7,
vscale=0.75,
headheight=13pt,
marginpar=2cm,
centering,
]{geometry}

\def\N{{\mathbb N}}
\def\Z{{\mathbb Z}}
\def\Q{{\mathbb Q}}

\def\C{{\mathbb C}}

\def\A{{\mathbb A}}

\def\cA{\mathcal{A}}
\def\cB{\mathcal{B}}

\def\cO{\mathcal{O}}
\def\cP{\mathcal{P}}
\def\cQ{\mathcal{Q}}
\def\cR{\mathcal{R}}
\def\cS{\mathcal{S}}

\def\fa{\mathfrak{a}}

\def\fm{\mathfrak{m}}
\def\fn{\mathfrak{n}}
\def\fp{\mathfrak{p}}
\def\fq{\mathfrak{q}}

\def\a{\alpha}

\def\Om{\Omega}

\def\.{\cdot}
\let\circum\^
\def\^{\widehat}
\def\~{\widetilde}

\def\ov{\overline}

\def\({\left(}
\def\){\right)}

\def\*{{}^*}

\newcommand{\cotimes}[1][]{\: \widehat{\otimes}_{#1} \:} 

\renewcommand{\and}{ \ \ \text{ and } \ \ }

\def\Jac{\mathrm{Jac}}

\def\MJ{\mathrm{MJ}}

\def\isom{\simeq} 

\DeclareMathOperator{\rk} {rk}

\DeclareMathOperator{\Spec} {Spec}
\DeclareMathOperator{\Spf} {Spf}

\DeclareMathOperator{\Sing} {Sing}

\DeclareMathOperator{\ord} {ord}

\DeclareMathOperator{\Sym} {Sym}

\DeclareMathOperator{\id} {id}

\DeclareMathOperator{\Hom} {Hom}

\DeclareMathOperator{\Cont} {Cont}
\DeclareMathOperator{\gr} {gr}

\DeclareMathOperator{\height} {ht}

\DeclareMathOperator{\trdeg}{trdeg}

\DeclareMathOperator*{\colim}{colim}
\DeclareMathOperator{\op}{op}
\DeclareMathOperator{\Pro}{Pro}

\DeclareMathOperator{\Set}{Set}

\def\embdim{\mathrm{edim}}
\def\embcodim{\mathrm{ecodim}}

\theoremstyle{plain}

\newtheorem{introtheorem}{Theorem}

\crefname{introtheorem}{Theorem}{Theorems}

\newtheorem{theorem}{Theorem}[section]
\newtheorem*{theorem*}{Theorem}
\newtheorem{proposition}[theorem]{Proposition}
\newtheorem{lemma}[theorem]{Lemma}
\newtheorem{corollary}[theorem]{Corollary}

\theoremstyle{definition}

\newtheorem{definition}[theorem]{Definition}

\theoremstyle{remark}

\newtheorem{remark}[theorem]{Remark}
\newtheorem{example}[theorem]{Example}
\newtheorem{question}[theorem]{Question}

\numberwithin{figure}{section}

\numberwithin{equation}{section}

\usepackage[
backend=biber,
style=numeric,
giveninits=true
]{biblatex}

\DefineBibliographyStrings{english}{%
	backrefpage = {cited on page},
	backrefpages = {cited on pages},
}

\DeclareSourcemap{
  \maps[datatype=bibtex]{
    \map{
        \step[fieldset=doi,null]
        \step[fieldset=url,null]
        \step[fieldset=isbn,null]
        \step[fieldset=issn,null]        
        \step[fieldset=note,null]
        \step[fieldset=pagetotal,null]
    }
  }
}

\newcommand{\SPC}[1]{\cite[\href{https://stacks.math.columbia.edu/tag/#1}{Tag #1}]{stacks-project}
}

\DeclareMathOperator{\Ad}{Ad}
\DeclareMathOperator{\Def}{Def}
\DeclareMathOperator{\Cpt}{Cpt}
\DeclareMathOperator{\Test}{Test}

\DeclareMathOperator{\TopRing}{TopRing}

\begin{document}

\title{The Drinfeld--Grinberg--Kazhdan theorem and embedding codimension of the arc space}
\author{Christopher Heng Chiu}
\address{University of Bern, Mathematical Institute, Sidlerstrasse 5, 3012 Bern, Switzerland.}
\email{christopher.chiu@unibe.ch}
\thanks{The author was supported by SNSF postdoc fellowship 217058.}

\begin{abstract}
    We prove an extension of the theorem of Drinfeld, Grinberg and Kazhdan to arcs with arbitrary residue field. As an application we show that the embedding codimension is generically constant on each irreducible subset of the arc space which is not contained in the singular locus. In the case of maximal divisorial sets, this relates the corresponding finite formal models with invariants of singularities of the underlying variety. We also prove an extension of a theorem by Bourqui and Sebag characterizing arcs of embedding codimension 0. 
\end{abstract}


\maketitle

\section*{Introduction}

Since the work of Nash, the geometry of the arc space $X_\infty$ of an algebraic variety $X$ has been known to encode information about the singularities of $X$. A key player in this connection is the notion of a \emph{maximal divisorial set}: to each divisorial valuation on $X$ one associates the closure $C_\nu(X)$ of the subset of arcs whose induced valuation agrees with $\nu$. For example, the Nash problem is characterizing those maximal divisorial sets which appear as irreducible components of the locus of arcs centered at $\Sing X$. In the context of birational geometry, it was proven in \cite{ELM04} that, for a smooth variety $X$, the codimension of a maximal divisorial subset $C_\nu(X)$ computes the discrepancy of $X$ along the divisorial valuation $\nu$. For a generalization to singular varieties see \cite{dFEI08}.

One of the starting points for this paper is the following, first proved over fields of characteristic $0$ in \cite{MR18} and then over perfect fields in \cite{dFD20,CdFD24}. If $\alpha_\nu$ denotes the generic point of $C_\nu(X)$, then one has
\begin{equation}
    \label{eq:intro1}
    \embdim(\cO_{X_\infty,\alpha_\nu}) = \^a_\nu(X)
\end{equation}
and
\begin{equation}
    \label{eq:intro2}
    \dim(\^{\cO_{X_\infty,\alpha_\nu}}) \geq a^{\MJ}_\nu(X),
\end{equation}
where $\^a_\nu(X)$ and $a^{\MJ}_\nu(X)$ are variants of discrepancies called \emph{Mather} and \emph{Mather-Jacobian (log) discrepancies}. This establishes a direct relation between invariants of the local ring at $\alpha_\nu$ and invariants of singularities of the base variety $X$. 

On the other hand, the theorem of Drinfeld, Grinberg and Kazhdan \cite{Dri02} says that for each $k$-rational $\alpha \in X_\infty$ not contained in $\Sing X$ one has a decomposition
\[
    \^{\cO_{X_\infty,\alpha}} \isom \^{\cO_{Z,z}} \cotimes k[[t_n \mid n \in \Z_{\geq0}]],
\]
where $Z$ is a scheme of finite type over $k$.
Any formal scheme isomorphic to the formal neighborhood of $Z$ at $z$ - that is, the formal spectrum of $\^{\cO_{Z,z}}$ - is called a finite formal model for $\alpha$. 
The question of how finite formal models are related to the singularities of $X$ is still wide open. One particular approach to providing an answer to this question is found in \cite{BMCS22} for toric varieties and in \cite{BMC23} for curves. In both papers, it was proved that for certain divisorial valuations $\nu$ on the respective variety $X$, there exists a nonempty open subset $U$ of $C_\nu(X)$ such for all $k$-rational $\alpha, \alpha' \in U$ their respective formal neighborhoods $(X_\infty,\alpha)$ and $(X_\infty,\alpha')$ are isomorphic. Furthermore, after changing the coefficient field, the formal scheme $(X_\infty,\alpha)$ is isomorphic to a formally smooth extension of $(X_\infty, \alpha_\nu)$, the formal neighborhood of the generic arc $\alpha_\nu$. Proving the first assertion is relatively straightforward; the hard part in both cases is obtaining an  explicit isomorphism for the second assertion. A major technical obstacle is that the results of \cite{BMCS22,BMC23} require a specific construction of a suitable coefficient field for $(X_\infty,\alpha_\nu)$, and it is unclear how to achieve this for a general variety $X$.

The strategy of this paper is different: we start by generalizing the statement of the Drinfeld--Grinberg--Kazhdan theorem to arcs with arbitrary residue field.

\begin{introtheorem}
    \label{intro:dgk}
    Let $X$ be a scheme locally of finite type over a perfect field $k$ and $\beta \in X_\infty \setminus (\Sing X)_\infty$. Then there exists a locally closed subscheme $V$ of $X_\infty$ containing $\beta$, a scheme $Z$ of finite type over $k$ and a morphism
    \[
       \mu \colon  V \to Z \times \A^\N
    \]
    such that for each $\alpha \in V$ there exists, up to finite separable extension of coefficient fields, an isomorphism between the formal neighborhoods $(X_\infty,\alpha)$ and $(Z \times \A^\N,\mu(\alpha))$.
\end{introtheorem}

We note here that in case $\alpha$ is $k$-rational, there is no nontrivial extension of coefficient fields needed - see \cref{r:dgk-k-rational} - and thus \cref{intro:dgk} does indeed fully recover the Drinfeld--Grinberg--Kazhdan theorem.
The proof of \cref{intro:dgk} takes up \cref{s:drinfeld} of this paper. We closely follow Drinfeld's proof in \cite{Dri02} by first constructing the \emph{scheme of formal models} $Z$ and then proving a bijection on the level of deformations. The final step involves results of \cite{CdFD24} on residue field extensions between arc spaces. In that way, we avoid explicitly constructing a coefficient field: we show that any choice of coefficient field for $\mu(\alpha)$ uniquely determines one for $\alpha$. One should compare the statement of \cref{intro:dgk} to attempts to extend the Drinfeld--Grinberg--Kazhdan theorem beyond the formal neighborhood, such as \cites{BNS16,Bou20,HW21}.

One consequence of \cref{intro:dgk} is that certain properties of $C_\nu(X)$ can be deduced from those of scheme of finite type over $k$. We will demonstrate this in the case of the \emph{embedding codimension}, which was introduced in \cite{CdFD22}. For arbitrary local rings $(R,\fm,K)$ it is defined as
\[
    \embcodim(R) \coloneqq \height(\ker(\Sym_K(\fm/\fm^2) \to \gr(R))).
\]
The name comes from the fact that for $R$ Noetherian one has $\embcodim(R) = \embdim(R) - \dim(R)$. One of the main results in \cite{CdFD22} was an explicit bound for the embedding codimension of arcs not contained in $\Sing X$, which together with the equation  \eqref{eq:intro1} immediately implies the inequality \eqref{eq:intro2}.

\begin{introtheorem}
    \label{intro:embcodim}
    Let $X$ be a variety over a perfect field $k$. If $W \subset X_\infty$ is any irreducible closed subset not contained in $(\Sing X)_\infty$, then the function
    \[
        \alpha \mapsto \embcodim( \cO_{X_\infty,\alpha}) \in \Z_{\geq0} \cup \{\infty\}
    \]
    is finite constant on a nonempty open subset of $W$.
\end{introtheorem}

In particular, \cite[Corollary 11.5]{dFD20} says that any maximal divisorial set is not contained in $(\Sing X)_\infty$ and thus \cref{intro:embcodim} holds in this case. That is, for any general $k$-rational arc in $C_\nu(X)$ the embedding codimension of any finite formal model equals that of the local ring at the maximal divisorial arc. This can be viewed as a step in understanding the precise relation of finite formal models of general arcs in $C_\nu(X)$ and invariants of $\nu$. One may reasonably expect that \cref{intro:dgk} allows to prove further results in this direction, for example directly expressing the Mather discrepancy in terms of finite formal models. An obvious obstruction is that finite formal models are not unique and neither is their (embedding) dimension. However, even when taking this into account, we will discuss in \cref{ss:max-div} how a precise relation is far from obvious. In particular, to our knowledge the question whether the formal neighborhood of a general $k$-rational arc in $C_\nu(X)$ is invariant of the choice of $\alpha$ is still open.

Recall that the embedding codimension of a Noetherian local ring is zero if and only if the ring is regular. Over perfect fields this can be generalized to non-Noetherian rings by replacing regularity with a suitable version of formal smoothness, see \cref{t:formal-smoothness-embcodim} and \cite[0, (19.5.4)]{egaiv1}. In keeping with the general philosophy of this paper we use the relation of embedding codimension at $k$-rational arcs and at more general points to extend a result by Bourqui and Sebag \cite{BS17b} as follows. The main theorem in \cite{BS17b} states that the local ring at a $k$-rational arc has embedding codimension $0$ if and only if the image of that arc is contained in a unique smooth formal branch. Using \cref{intro:embcodim} we can prove:

\begin{introtheorem}
    \label{intro:smooth-arcs}
    Let $k$ be a perfect field, $X$ a variety over $k$ and $\alpha \in X_\infty$. Assume that either
    \begin{enumerate}
        \item $k$ is uncountable, or
        \item $\alpha$ is \emph{stable} (see \cref{d:stable}).
    \end{enumerate}
    If $X$ is geometrically unibranch at $\alpha(0)$, then $\embcodim(\cO_{X_\infty,\alpha}) = 0$ if and only if $X$ is smooth at $\alpha(0)$.
\end{introtheorem}

The assumption on the field $k$ (or that $\alpha$ is stable) is needed in order to ensure the existence of enough rational points. In the case where $\alpha$ is stable \cref{intro:smooth-arcs} addresses a variant of a question by Reguera in \cite[Question 2.10]{Reg18}. Reguera's original question assumed that $X$ is unibranch at $\alpha(0)$ (a weaker assumption) and $\cO_{X_\infty,\alpha}$ regular (a stronger assumption). We explain in \cref{p:unibranch-degen} how \cref{intro:smooth-arcs} does not in fact fully resolve \cite[Question 2.10]{Reg18}.

\subsection*{Acknowledgments}

We are very grateful to David Bourqui and Mario Mor\'an Can\'on for pointing out that the main result of the paper can be used to address \cite[Question 2.10]{Reg18}, which lead to \cref{intro:smooth-arcs}. We thank Tommaso de Fernex and Roi Docampo for the useful discussions and in particular for pointing out the observation in \cref{r:reduced-main-cpt}. Finally, we thank the referee for their careful reading of the paper and their many suggestions which helped improve the exposition of the paper significantly.

\section{Preliminaries}
\label{s:preliminaries}

Throughout this paper we fix $k$ to be a perfect field. By variety over $k$ we mean an integral separated scheme of finite type over $k$. If $f \colon X \to Y$ is a morphism of affine schemes, then the corresponding map between coordinate rings will be denoted by $f^\sharp$.

\subsection{Weierstrass preparation and division theorem}
\label{ss:weierstrass}

Two of the main algebraic ingredients in the proof of \cref{intro:dgk} are the Weierstrass preparation and division theorems. As in this paper we go beyond \cite{Dri02} we include the precise statements used later in \cref{s:drinfeld}.

\begin{definition}
    Let $(A,\fm)$ a local ring. We call a polynomial of the form
    \[
        q(t) = t^d + q_{d-1} t^{d-1} + \ldots + q_1 t + q_0
    \]
    with $q_j \in \fm$, $j=0,\ldots,d-1$, a \emph{Weierstrass polynomial} of degree $d$.
\end{definition}

\begin{proposition}[Weierstrass preparation and division]
    Let $(A,\fm)$ be a local complete ring with residue field $K$.
    \begin{enumerate}
        \item Let $f(t) \in A[[t]]$ and write $f_0 \in K[[t]]$ for its image modulo $\fm$. If $\ord_t f_0(t) = d < \infty$, then there exists a Weierstrass polynomial $q(t)$ and a unit $u(t) \in A[[t]]^*$ with $f(t) = u(t) q(t)$.
        \item Let $q(t)$ be a Weierstrass polynomial of degree $d$. For every $f(t) \in A[[t]]$ there exists unique $g(t) \in A[[t]]$ and $r(t) \in A[t]_{<d}$ with $f(t) = g(t) q(t) + r(t)$.
    \end{enumerate}
\end{proposition}

For a proof we refer the reader to \cite[VII, \S3, 8-9]{Bou}. A crucial property of Weierstrass polynomials with coefficients in a complete local ring is the following.

\begin{lemma}
    \label{l:weierstrass-polynomial}
    Let $(A,\fm)$ be a local ring and $q(t) \in A[t]$ a Weierstrass polynomial of degree $d$. Assume that $A$ is separated, i.e.\ $\bigcap_n \fm^n = (0)$.
    \begin{enumerate}
        \item \label{eq1:weierstrass-polynomial} We have $q A[[t]] \cap A[t] = q A[t]$. 
        \item \label{eq2:weierstrass-polynomial} The element $q$ is a regular element in $A[t]$ and $A[[t]]$. In particular, every element $f(t) \in A[[t]]$ with $f_0(t) \in K[[t]]$ nonzero is a regular element.
    \end{enumerate}
    If in addition $(A,\fm)$ is a complete local ring, then every $f(t) \in A[[t]]$ whose image $f_0(t) \in K[[t]]$ is nonzero is a regular element of $A[[t]]$.
\end{lemma}

\begin{proof}
    For \eqref{eq1:weierstrass-polynomial}, assume that $p(t) = q(t) u(t) \in A[t]$ with $u(t) \in A[[t]]$. If $e = \deg_t p(t)$, then for all $i>e$ we have
    \[
        0 = u_{i-d} + q_{d-1}u_{i+1-d} + \ldots + q_1 u_{i-1} + q_0 u_i.
    \]
    Since $q_{d-1},\ldots,q_0 \in \fm$, it follows that for $j>e$ we have $u_{j-d} \in \fm$ and thus $u_{j-d} \in \fm^n$ for all $n$. 

    \eqref{eq2:weierstrass-polynomial} is clear for $q \in A[t]$. For $q \in A[[t]]$, the argument is almost identical to the one above. Finally, the last assertion follows from \eqref{eq2:weierstrass-polynomial} and Weierstrass preparation.
\end{proof}

\begin{remark}
    If $X$ is an algebraic variety and $\alpha \in X_\infty$ an arc such that $\alpha(0) \in X$ is singular, then the local ring at $\alpha$ is in general not separated. This is demonstrated in \cite[Example 5.13 and Proposition 8.3]{CdFD24} for $X$ the cuspidal plane curve $x^3 - y^2 = 0$. We will see another example in \cref{p:unibranch-degen}.
\end{remark}

\subsection{Describing completions via deformations}
\label{ss:formal-schemes}

A key observation in Drinfeld's proof in \cite{Dri02} is the following. If $(R,\fm)$ is a local $k$-algebra with residue field $k$ its $\fm$-adic completion $\^R$ is an inverse limit of \emph{test rings}; that is, local $k$-algebras with residue field $k$ and nilpotent maximal ideal. Then the ring $\^R$ is determined by the functor of points of $A$ restricted to test rings. For our proof of \cref{intro:dgk} we need to both allow arbitrary residue fields and rings which do not necessarily arise as completions of local rings. As such, it makes sense to instead think of the underlying inverse system as a pro-object on a category of test rings with a fixed identification of their residue field.

\begin{definition}
    Let $K/k$ be a field extension. The category $\Test_K$ of $K$-test rings has as objects local $k$-algebras $(A,\fm_A)$ together with a local $k$-surjection $\sigma_A \colon A \to K$ such that $\fm_A^n = 0$ for some $n\geq 0$. A morphism $A \to B$ of $K$-test rings is a local $k$-algebra map $\varphi \colon A \to B$ such that the diagram
    \[
        \begin{tikzcd}
            A \ar[r, " \varphi"] \ar[d, "\sigma_A"] & B \ar[d, "\sigma_B"] \\
            K \ar[r, "\id_K"] & K
        \end{tikzcd}
    \]
    commutes.
\end{definition}

We will often omit the surjection $\sigma_A$ from the notation if it is clear from the context.

\begin{remark}
    Any $(A,\sigma_A) \in \Test_K$ is a complete local ring over a field $k$. In particular, there exists a coefficient field $\iota_A \colon K \to A$. If $j \colon K \to L$ is a field extension, then for any choice of coefficient field for $A$ the base change $A \otimes_K L$ (equipped with the natural map) gives an element of $\Test_L$.
\end{remark}

We denote by $\N$ the category whose objects are elements $i\in \N$ and with exactly one morphism $i \to j$ if $i \leq j$.

\begin{definition}
    Let $\Pro(\Test_K)$ be the category of pro-objects of $\Test_K$. We define $\Cpt_K$ to be the full subcategory of $\Pro(\Test_K)$ consisting of \emph{surjective} pro-objects of $\Test_K$ which are indexed by $\N$. That is, objects of $\Cpt_K$ are functors $\cA \colon \N^{\op} \to \Test_K$, with components $A_i \coloneqq \cA(i)$, such that $A_j \to A_i$ is surjective for $i\leq j$. For the morphisms, we have
    \[
        \Hom_{\Cpt_K}(\cA,\cB) = \varprojlim_j \varinjlim_i \Hom_{\Test_K}(A_i,B_j).
    \]
\end{definition}

We call a functor $\Test_K \to \Set$ pro-representable if it is isomorphic to $\colim_i h_{A_i}$, where $h_{A_i} \coloneqq \Hom_{\Test_K}(A_i,-)$ and $\cA = (A_i)_{i\in \N} \in \Pro(\Test_K)$. By the Yoneda lemma, the contravariant functor $\cA \mapsto \varinjlim_i h_{A_i}$ defines an equivalence between the category of pro-objects and the full subcategory of pro-representable functors $\Test_K \to \Set$ (see \cite{fgaiii}).

Precomposition with the functor $\{\star\} \to \N$ gives a fully faithful functor $\Test_K \to \Cpt_K$, by which we can consider each $A \in \Test_K$ as an object in $\Cpt_K$ (a constant inverse system). Clearly $K$ is a final object in $\Cpt_K$.

\begin{lemma}
    Let $\cA \in \Cpt_K$, then there exists a section $\iota_{\cA} \colon K \to \cA$.
\end{lemma}

We call such a section $\iota_{\cA}$ a \emph{coefficient field} of $\cA$.

\begin{proof}
    Consider the diagram
    \[
        \begin{tikzcd}
            k \ar[r] \ar[d] & A_i \ar[d] \\
            K \ar[r] & A_{i-1}.
        \end{tikzcd}
    \]
    Note that $A_i \to A_{i-1}$ is surjective and its kernel $I \subset A_i$ satisfies $I^n = 0$ for some $n\geq 0$. Since $k$ is perfect we have that $k \to K$ is formally smooth and the result follows.
\end{proof}

In other words, a choice of coefficient field $\iota_{\cA}$ for $\cA$ corresponds to compatible choices of coefficient fields $\iota_{A_i}$ for $A_i$.

Denote by $\TopRing_k$ the category of topological rings over $k$ (considered with the discrete topology). For $\cA \in \Cpt_K$ the inverse limit $\varprojlim_i A_i$ endowed with the induced topology is an object of $\TopRing_k$. Moreover, for $B \in \Test_K$ we have a natural map
\[
    \varinjlim_i \Hom_{\Test_K}(A_i,B) \to \Hom_{\TopRing_k}(\varprojlim_i A_i,B),
\]
which is injective as $\cA$ has surjective transition maps. Not that the above map is not surjective: a continuous map $A_i \to B$ might not be compatible with the choice of identifications $\sigma_{A_i}$ and $\sigma_B$. To summarize, we have:

\begin{lemma}
    \label{l:limit-top-rings}
    The assignment $\cA \mapsto \varprojlim_i A_i$ gives a faithful (but not full) functor $\Cpt_K \to \TopRing_k$. Any choice of coefficient field of $\cA$ gives a coefficient field for $\varprojlim_i A_i$.
\end{lemma}

For a choice of coefficient field for $A \in \Test_K$ and a field extension $j \colon K \to L$ the base change $A \otimes_K L$ has a canonical surjection $\sigma_{A \otimes_K L} = (j \circ \sigma_A) \otimes \id_L$ and thus becomes an object of $\Test_L$. More generally, we have the following.

\begin{lemma}
    \label{l:base-change-cpt}
    Let $\cA = (A_i)_{i\in\N} \in \Cpt_K$ and fix a choice of coefficient field for $\cA$. Let $j \colon K \to L$ be a field extension. Define $\cA \otimes_K L$ by
    \[
        (\cA \otimes_K L)(i) \coloneqq A_i \otimes_K L.
    \]
    Then $\cA \otimes_K L \in \Cpt_L$. Moreover
    \[
        \varprojlim_i (A_i \otimes_K L) \isom (\varprojlim_i A_i) \^\otimes_K L,
    \]
    where on the right hand side we consider the completed tensor product as cofibered coproduct in $\TopRing_k$ (see \cite[0, \S 7.7]{egai}).
\end{lemma}

\begin{proof}
    We have that $\cA \otimes_K L \in \Cpt_L$ as taking tensor products is right exact. The second assertion is obvious from the definition of the completed tensor product.
\end{proof}

Let $S$ be any $k$-algebra together with a $k$-homomorphism $\tau \colon S \to L$ where $L$ is a field. For $(A,\sigma_A) \in \Test_L$ an \emph{$A$-valued deformation} of the pair $(S,\tau)$ is a commutative diagram
\[
    \begin{tikzcd}
        S \ar[r, "\varphi"] \ar[d, "\tau"'] & A \ar[ld, "\sigma_A"] \\
        L. &
    \end{tikzcd}
\]
We denote the set of $A$-valued deformations of $(S,\tau)$ by $\Def_{(S,\tau)}(A)$. This gives a functor $\Def_{(S,\tau)} \colon \Test_L \to \Set$.

\begin{example}
    Let $S = k[x_1,\ldots,x_n]/(f_1,\ldots,f_r)$ and $\tau \colon S \to L$ a $k$-homomorphism with $L$ a field. Let $(A,\sigma_A) \in \Test_L$ and $\fm_A$ be the maximal ideal of $A$. We can identify $A$-valued deformations of $(S,\tau)$ with tuples $(a_1,\ldots,a_n) \in A^n$ such that $\sigma_A(a_i) = \tau(x_i)$ and $f_l(a_1,\ldots,a_n) = 0$ for all $l = 1,\ldots,r$.
\end{example}

We will mostly apply the above notions in the following case. Let $(R,\fm)$ be a local $k$-algebra with residue field isomorphic to $K$. Choosing an isomorphism $R/\fm \isom K$ defines a surjective $k$-algebra homomorphism $\sigma \colon R \to K$ and thus an object $\cR \in \Cpt_K$ via $\cR(i) \coloneqq R/\fm^{i+1}$, with $\sigma_{R_i} \colon R/\fm^i \to K$ the map induced by $\sigma$. The image of $\cR$ in $\TopRing_k$ is isomorphic to the $\fm$-adic completion $\widehat{R}$ of $R$.

Now fix a field extension $j \colon K \to L$ over $k$ and write $\sigma^{L} \coloneqq j \circ \sigma$. Let $(A,\sigma_A) \in \Test_L$ and take any $A$-valued deformation of the pair $(R,\sigma^L)$
\[
    \begin{tikzcd}
        R \ar[r, "\varphi"] \ar[d, "\sigma^L"'] & A \ar[ld, "\sigma_A"] \\
        L. &
    \end{tikzcd}
\]
Note that then $\varphi$ is automatically local and so in particular continuous for the $\fm$-adic topology on $R$ (and the discrete topology on $A$). If the maps $\sigma$ and $j$ are clear from the context, we will also write $\Def_{R,L/K} \coloneqq \Def_{(R,\sigma^L)}$. Moreover, in case $j = \id_K$ we will abbreviate to $\Def_{R}$ and call elements in $\Def_R(A)$ just $A$-valued deformations of $R$.

The following is a key observation in extending the proof of Drinfeld in \cite{Dri02} to arcs which are not $k$-rational (see \cref{ss:formal-neighborhood}).

\begin{proposition}
    \label{p:deformations-finite-separable}
    Let $(R,\fm)$ be a local ring with together with a local surjection $\sigma \colon R \to K$ over $k$. Let $j \colon K \to L$ a finite separable field extension and write $\sigma^L \coloneqq j \circ \sigma$. Let $\cR \in \Cpt_K$ corresponding to $R$ and choose a coefficient field $\iota \colon K \to \cR$. Then the natural map
    \begin{equation}
        \label{eq:deformations-finite-separable}
        \varinjlim_i \Hom_{\Test_L}(R_i \otimes_K L,-) \to \Def_{(R,\sigma^L)}
    \end{equation}
    is an isomorphism.
\end{proposition}

\begin{proof}
    Let $(A,\sigma_A) \in \Test_L$ and denote by $\fm_A$ the maximal ideal of $A$. We have to show that \eqref{eq:deformations-finite-separable} is bijective when evaluated at $A$. An element of the left hand side is given by a pair $(\psi_1,\psi_2)$ fitting into the diagram
    \[
        \begin{tikzcd}
            K \ar[r,"\iota_i"] \ar[d,"j"'] & R_i \ar[d,"\psi_1"] \ar[r,"\sigma_{R_i}"] & K \ar[d,"j"] \\
            L \ar[r,"\psi_2"] & A \ar[r, "\sigma_A"] & L
        \end{tikzcd}
    \]
    and its image in $\Def_{(R,\sigma^L)}$ is given by the composition $\varphi \coloneqq \pi_i \circ \psi_1$, where $\pi_i \colon R \to R_i$. Since $\pi_i$ is surjective we have that $\pi_i \circ \psi_1 = \pi_i \circ \psi'_1$ implies $\psi_1 = \psi'_1$. Conversely, since $\fm_A$ is nilpotent every $\varphi \colon R \to A$ in $\Def_{(R,\sigma^L)}$ factors through some $\pi_{i'} \colon R \to R_{i'}$. It is thus sufficient to show that for any $\psi_1$ with $\sigma_A \circ \psi_1 = j \circ \sigma_{R_i}$, there exists a unique coefficient field $\psi_2$ which makes the above diagram commute.

    To that avail, any such coefficient field corresponds to a diagonal arrow in the diagram
    \[
        \begin{tikzcd}
            K \ar[r,"\psi_1 \circ \iota_i"] \ar[d, "j"'] \ar[rd, dashed] & A \ar[d,"\sigma_A"] \\
            L \ar[r,"\id_L"] & L.
        \end{tikzcd}
    \]
    Since $j \colon K \to L$ is finite separable it is formally \'etale and there exists a unique diagonal arrow in the diagram.
\end{proof}

\subsection{The arc space of an algebraic variety}

We will briefly introduce some elementary facts from the theory of arc spaces and fix some notation. For a more comprehensive treatment we refer the reader to the various introductory texts available in the literature.

Let $X$ be any scheme over $k$. The arc space $X_\infty$ of $X$ is obtained as the limit $X_\infty = \varprojlim_n X_n$, with the $n$-th jet space $X_n$ defined via
\[
    \Hom_k(Z,X_n) \isom \Hom_k(Z \times_k \Spec k[t]/(t^{n+1}),X).
\]
Note that if $X$ is affine, then so is $X_\infty$. In fact, writing $X = \Spec R$, we have that $X_\infty = \Spec R_\infty$, where $R_\infty$ denotes the algebra of higher derivations \cite{Voj07}. Explicitly, one has a presentation of $R_\infty$ of the form
\[
    R [x^{(i)} \mid x \in R, i\in \Z_{\geq 1 }] \twoheadrightarrow R_\infty,
\]
where we will write $x^{(0)} \in R_\infty$ for the image of $x \in R$. The algebra $R_\infty$ comes equipped with a universal higher derivation $R \to R_\infty[[t]]$, $x \mapsto \sum_{i \geq 0} x^{(i)} t^i$, and satisfies
\[
    \Hom_k(R_\infty, R) \isom \Hom_k(\Spec R, X_\infty) \isom \Hom_k(\Spec R[[t]], X).
\]
More generally, for an arbitrary scheme $X$ one can deduce from the affine case that $X_\infty$ is a $k$-scheme which satisfies
\[
    \Hom_k(\Spec K, X_\infty) \isom \Hom_k(\Spec K[[t]],X)
\]
for any field extension $K \supset k$. Thus we identify points $\alpha \in X_\infty$ with the corresponding morphism $\alpha \colon \Spec k_\alpha[[t]] \to X$. Writing $\Spec K[[t]] = \{0,\eta\}$, the projection $\pi \colon X_\infty \to X$ is given by $\alpha \mapsto \alpha(0)$. Moreover, for any closed subset $Z \subset X$ we have $\alpha \in Z_\infty$ if and only if $\alpha(\eta) \in Z$.

Now let $X$ be a variety. Any arc $\alpha \in X_\infty$ defines a semi-valuation $\ord_\alpha$ on $\cO_X(U)$ with $\alpha(0) \in U$ affine via $\ord_\alpha(f) \coloneqq \ord_t(\alpha^\sharp(f))$. If $\alpha(\eta)$ is the generic point of $X$, then $\ord_\alpha$ actually gives a $\Z$-valued valuation of the function field $k(X)$ of $X$. Note that if $\alpha'$ specializes to $\alpha$, then $\alpha'(0)$ specializes to $\alpha(0)$ and for every affine $U \subset X$ with $\alpha(0)\in U$ we have $\ord_{\alpha'} \leq \ord_\alpha$ on $\cO_X(U)$.

\section{On the Drinfeld--Grinberg--Kazhdan theorem}
\label{s:drinfeld}

In this section we revisit Drinfeld's proof in \cite{Dri02} and extend it to prove \cref{intro:dgk}. The core of the argument follows the one in \cite{Dri02} for $k$-rational arcs: the bijection between deformations constructed in \cref{ss:bijection-deformations}. We give a full proof here to demonstrate the validity of the argument when passing to general points. To deduce the isomorphism of formal neighborhoods we make use of results from \cite{CdFD24} on residue field extensions at the level of arc spaces.

\subsection{Reduction to the case of complete intersection}
\label{ss:complete-intersection}

To prove \cref{intro:dgk}, we may assume that $X$ is affine. The next step is a standard argument to reduce to the case of a complete intersection. This was used in \cite{Dri02} and explained in more detail in \cite[Section 4.2]{BS17}. For the reader's convenience we recall the proof here to show that it extends from a single arc to an open neighborhood of $X_\infty$. We first recall the following nonstandard notation from the introduction.

\begin{definition}
    Let $X$ be a scheme over $k$ and $x\in X$. Then we write
    \[
        (X,x) \coloneqq \Spf \widehat{\cO_{X,x}}
    \]
    and call it the \emph{formal neighborhood} of $X$ at $x$.
\end{definition}

Note that, if $x$ is a closed point of $X$, then $(X,x)$ is isomorphic to the formal completion of $X$ along $x$.

\begin{proposition}
    \label{p:complete-intersection}
    Let $X$ be an affine scheme of finite type over $k$ and $\beta \in X_\infty \setminus (\Sing X)_\infty$. Then there exists a closed immersion $X \to X'$ with
    \[
        X' = \Spec k[x_1,\ldots,x_N]/(f_1,\ldots,f_r),
    \]
    an $r \times r$-minor $\delta$ of the Jacobian $(\frac{\partial f_i}{\partial x_j})_{i,j}$ and an open neighborhood $U \subset X_\infty$ of $\beta$ satisfying the following. For each $\alpha \in U$, we have $\alpha^\sharp(\delta) \neq 0$ and the induced map of formal neighborhoods
    \[
        (X_\infty,\alpha) \to (X'_\infty,\alpha)
    \]
    is an isomorphism.
\end{proposition}

It follows from \cref{ss:formal-schemes} that the local ring $\cO_{X_\infty,\alpha}$ defines an object of $\Cpt_K$, where $K \coloneqq k_\alpha$, and similar for $\cO_{X'_\infty,\alpha}$. We note that the isomorphism of \cref{p:complete-intersection} comes from an isomorphism in $\Cpt_K$; that is, it is compatible with the obvious identification of residue fields.

\begin{proof}
    Let $X = \Spec k[x_1,\ldots,x_N]/\fa$, then the ideal of $\Sing X$ is the radical of the ideal generated by elements of the form $h \delta$, where $\delta$ is a minor of the Jacobian matrix of some elements $f_1,\ldots,f_r \in \fa$ and $h \in ((f_1,\ldots,f_r) : \fa)$. By assumption $\beta$ is not contained in $\Sing X$ and hence there exist $h,f_1,\ldots,f_r$ and $\delta$ as before such that $\beta^\sharp(h\delta) \neq 0$.
    Let $X' \coloneqq \Spec k[x_1,\ldots,x_N]/(f_1,\ldots,f_r)$, which contains $X$ as a closed subscheme. Consider the open subset
    \[
        X^{h,\delta} \coloneqq \{\alpha \in X \colon \alpha^\sharp(h\delta)\neq0 \}
    \]
    of $X$. We claim that for all $\alpha \in X^{h,\delta}$ the natural map of formal neighborhoods
    \[
        (X_\infty,\alpha) \to (X'_\infty,\alpha)
    \]
    is an isomorphism. Indeed, since this map induces an isomorphism of residue fields, by \cref{l:limit-top-rings} it suffices to show that every $A$-valued deformation of $\cO_{X'_\infty,\alpha}$ lifts via the above morphism to an $A$-valued deformation of $\cO_{X_\infty,\alpha}$, where $(A,\fm_A)$ is a $k_\alpha$-test ring. So let $\widetilde\alpha \colon \Spec A[[t]] \to X'$ be given by $x_1(t),\ldots,x_N(t) \in A[[t]]$ satisfying $f_i(x_1(t),\ldots,x_N(t)) = 0$ for $i=1,\ldots,r$ and $x_j(t) \equiv x_j^0(t)$ mod $\fm_A$, where $x_1^0(t),\ldots,x_N^0(t) \in k_\alpha[[t]]$ are the images of the $x_j$'s under $\alpha^\sharp$. In particular, since
    \[
        h(x^0(t)) \delta(x^0(t)) \neq 0,
    \]
    by \cref{l:weierstrass-polynomial} the element $h(x(t)) \in A[[t]]$ is regular. By definition that implies that $f(x(t)) = 0$ for all $f \in \fa$ and hence $\widetilde{\alpha}$ lifts to $X$.
\end{proof}

From now on we will fix the following situation. Write $x = (x_1,\ldots,x_n)$, $y = (y_1,\ldots,y_m)$. We assume $X = V(f_1,\ldots,f_m) \subset \A^{n+m}$ with $f=(f_1,\ldots,f_m) \in k[x,y]^m$. Let $Df= (\frac{\partial f_i}{\partial y_j})_{i,j \leq m}$ and $\delta = \det Df$. For $d \geq 0$ define
\[
    X_\infty^{\delta,d} \coloneqq \{\alpha \in X_\infty \mid \ord_\alpha \delta = d\}.
\]
Note that $X_\infty^{\delta,d}$ is a locally closed subset of $X_\infty$, given by
\[
    X_\infty^{\delta,d} = V(\delta^{(0)},\ldots,\delta^{(d-1)}) \cap D(\delta^{(d)}).
\]
For any matrix $M$ with coefficients in $R$ write $\Ad (M)$ for the adjoint matrix of $M$. Note that the coefficients of $\Ad(M)$ are polynomials in the coefficients of $M$.

\subsection{The scheme of formal models}
\label{ss:formal-model}

Let $\cQ_d$ be the scheme of monic polynomials of degree $d$ in one variable $t$. In other words, for each $k$-algebra $R$ we have
\[
    \cQ_d(R) = \{q = t^d + q_{d-1}t^{d-1} + \ldots + q_1 t + q_0 \mid q_i \in R\}.
\]
Clearly $\cQ_d \isom \A^d$. Similarly, write $\cP_{<e}$ for the scheme of polynomials in one variable of degree $<e$; that is, $\cP_{<e}(R) = R[t]_{<e}$ for every $k$-algebra $R$. Finally write $\cP_\infty = (\A^1)_\infty$, so $\cP_\infty(R) = R[[t]]$. 

\begin{lemma}
    \label{l:formal-model-scheme}
    Consider the scheme $W = \cQ_d \times \cP_{<2d}^n \times \cP_{<d}^m$. If $R$ is a $k$-algebra, then $Z(R)$ is the set
    \[
        \{(q,\ov{x},\ov{y}) \in W(R) \mid \delta(\ov{x},\ov{y}) \in q R[t], \Ad(Df(\ov{x},\ov{y})) \cdot f(\ov{x},\ov{y}) \in q^2 R[t]^m \}.
    \]
    The functor $R\mapsto Z(R)$ is representable by a scheme $Z$ of finite type over $k$.
\end{lemma}

\begin{proof}
    The idea is to interpret $\ov{x}, \ov{y}$ as residue classes modulo $q$. That is, for every polynomial $p \in R[t]$ there exists a unique polynomial $r \in R[t]_{<d}$ such that $p \equiv r \mod q R[t]$. Thus we can write
    \[
        \delta(\ov{x},\ov{y}) \equiv r_\delta \mod q R[t]
    \]
    with $r_\delta \in R[t]_{<d}$. Observe that the coefficients $r_{\delta,i}$, $i<d$, of $r_\delta$ are polynomials over $k$ in the coefficients of $\ov{x}, \ov{y}$ and $q$. Thus the $r_{\delta,i}$ define regular functions on $W$ and the system of equations $r_{\delta,i} = 0$, $i<d$, is equivalent to $\delta(\ov{x},\ov{y}) \in q R[t]$. Similarly arguing for the relation
    \[
        \Ad(Df(\ov{x},\ov{y})) \cdot f(\ov{x},\ov{y}) \in q^2 R[t]^m
    \]
    one obtains equations defining $Z$ as a subscheme of $W$.
\end{proof}

Now let $(x,y) \in X_\infty^{\delta,d}(R)$. We may write $\delta(x,y) = t^d u$ with $u \in R[[t]]^*$. Write  $x = \ov{x} + t^{2d} \xi$, $y = \ov{y} + t^d \theta$ with $\ov{x} \in \cP_{<2d}^n(R)$, $\xi \in \cP_\infty^m(R)$ and similarly $\ov{y} \in \cP_{<d}^m(R)$, $\theta \in \cP_\infty^m(R)$. We claim
\[
    (x,y) \mapsto (t^d,\ov{x},\ov{y};\xi)
\]
defines a morphism $\mu \colon X_\infty^{\delta,d} \to Z \times \cP_\infty^n$. Indeed, we have
\[
    \delta(\ov{x},\ov{y}) \equiv 0 \mod t^d.
\]
Moreover, using Taylor expansion we have
\[
    0 = f(x,y) \equiv f(\ov{x},\ov{y}) + t^d Df(\ov{x},\ov{y}) \cdot \theta \mod t^{2d} R[t]^m
\]
and hence
\[
    \Ad(Df(\ov{x},\ov{y})) \cdot f(\ov{x},\ov{y}) \equiv t^d \delta(\ov{x},\ov{y}) \theta \equiv 0 \mod t^{2d} R[t]^m.
\]
Note that here we used crucially that $t^e R[[t]] \cap R[t] = t^e R[t]$.

\subsection{A bijection of deformations}
\label{ss:bijection-deformations}

Let $X$ be a scheme over $k$, $K$ a field extension of $k$ and $(A,\sigma_A) \in \Test_K$. Let $x_K$ be a $K$-point of $X$. Following \cref{ss:formal-schemes} we define an \emph{$A$-valued deformation of the pair $(X,x_K)$} to be a commutative diagram
\[
    \begin{tikzcd}
        \Spec A \ar[r, "\widetilde{x}"] & X \\
        \Spec K \ar[u] \ar[ru, "x_K"'] & 
    \end{tikzcd}
\]
with the vertical arrow induced by $\sigma_A$. We will write $\Def_{(X,x_K)} \colon \Test_K \to \Set$ for the corresponding functor. Let $x \in X$ be the image of $x_K$ and $k_x$ its residue field. As each $A$-valued deformation $\widetilde{x}$ factors through $\Spec \cO_{X,x}$ we have
\[
    \Def_{(X,x_K)} = \Def_{\cO_{X,x},K/k_x}.
\]
Now, if $\alpha_K$ is a $K$-point of $X_\infty$, then we will identify any $A$-valued deformation $\~\alpha$ of $(X_\infty,\alpha_K)$ with the corresponding diagram
\[
    \begin{tikzcd}
        \Spec A[[t]] \ar[r, "\~\alpha"] & X \\
        \Spec K[[t]] \ar[u] \ar[ru, "\alpha_K"'] & 
    \end{tikzcd} 
\]
where $\a_K$ here denotes the corresponding $K[[t]]$-point of $X$.

The following is the main content of Drinfeld's proof in \cite{Dri02}.

\begin{theorem}
    \label{t:deformations}
    Let $X$ be as in \cref{ss:complete-intersection} and $\alpha_K \colon \Spec K \to X_\infty^{\delta,d}$ a $K$-point of $X_\infty^{\delta,d}$ (and by slight abuse of notation of $X_\infty$). Let $\mu \colon X_\infty^{\delta,d} \to Z \times \cP_\infty^n$ be the morphism from \cref{ss:formal-model}. Write $\gamma_K$ for the $K$-point defined by $\mu \circ \alpha_K$. Then there exists a natural isomorphism
    \[
        \Phi \colon \Def_{(X_\infty,\alpha_K)} \to \Def_{(Z\times \cP_\infty^n,\gamma_K)}.
    \]
\end{theorem}

\begin{proof}
    Let $(A,\fm) \in \Test_K$ be a $K$-test ring. Note that the morphism $\alpha_K \colon \Spec K \to X_\infty$ corresponds to $\alpha(t) = (x_0(t),y_0(t)) \in K[[t]]^{n+m}$ with
    \[
        f(x_0(t),y_0(t)) = (f_1(x_0(t),y_0(t)),\ldots,f_m(x_0(t),y_0(t))) = 0
    \]
    and satisfying 
    \[
        \delta(x_0(t),y_0(t)) = \det \left(\frac{\partial f_i}{\partial x_j}(x_0(t),y_0(t))\right)_{i,j\leq m} = t^d u_0(t)
    \]
    with $u_0(t) \in K[[t]]$ invertible. Then any $A$-valued deformation $\widetilde{\alpha}$ of $(X_\infty,\alpha_K)$ is given by $\widetilde{\alpha}(t) = (x(t),y(t)) \in A[[t]]^{n+m}$ such that $x(t) \equiv x_0(t), y(t) \equiv y_0(t) \mod \fm$ and satisfying
    \begin{gather*}
        f(x(t),y(t)) = 0, \\
        \delta(x(t),y(t)) = q(t) u(t),
    \end{gather*}
    where $q(t) \in A[t]$ is a Weierstrass polynomial of degree $d$ and $u(t) \in A[[t]]$ is a unit. Note that the existence and uniqueness of $q(t), u(t)$ follows from the Weierstrass preparation theorem, and similary we get $u(t) \equiv u_0(t) \mod \fm$.

    By \cref{l:formal-model-scheme} and the definition of $\mu$ the morphism $\gamma_K \colon \Spec K \to Z \times \cP_\infty^n$ corresponds to $(t^d,\ov{x}_0(t),\ov{y}_0(t);\xi_0(t))$ with $\ov{x}_0(t) \in K[t]^n_{<2d}$, $\ov{y}_0(t)\in K[t]^m_{<d}$ and $\xi_0(t) \in K[[t]]^n$ such that
    \[
        x_0(t) = \ov{x}_0(t) + t^{2d} \xi_0(t), \: y_0(t) \equiv \ov{y}_0(t) \mod t^d K[[t]]^m.
    \]
    An $A$-valued deformation $\widetilde{\gamma}$ of $(Z\times\cP_\infty^n,\gamma_K)$ is then given by $(q(t),\ov{x}(t),\ov{y}(t),\xi(t))$ with $q(t) \in A[t]$ a Weierstrass polynomial of degree $d$, $\ov{x}(t) \in A[t]^n_{<2d}$, $\ov{y}(t) \in A[t]^m_{<d}$ and $\xi(t) \in A[[t]]^n$. These satisfy the following conditions: first,
    \[
        \ov{x}(t) \equiv \ov{x}_0(t), \ov{y}(t) \equiv \ov{y}_0(t), \xi(t) \equiv \xi_0(t) \mod \fm,
    \]
    and
    \begin{gather*}
        \delta(\ov{x}(t),\ov{y}(t))  \in q A[t], \\
        \Ad(Df(\ov{x}(t),\ov{y}(t))) \cdot f(\ov{x}(t),\ov{y}(t))  \in q^2 A[t]^m.
    \end{gather*}
    Now we define the map $\Phi_A \colon \Def_{(X_\infty,\alpha_K)}(A) \to \Def_{(Z \times \cP_\infty^n,\gamma_K)}(A)$ as follows. Given $\widetilde{\alpha}(t) = (x(t),y(t)) \in A[[t]]^{n+m}$, we let $\delta(x(t),y(t)) = q(t) u(t)$ for $q(t) \in A[t]$ a Weierstrass polynomial of degree $d$ and $u(t) \in A[[t]]$ invertible. Then write $x(t) = \ov{x}(t) + q^2(t) \xi(t)$, $y(t) = \ov{y}(t) + q(t) \theta(t)$ for $\xi(t), \theta(t) \in A[[t]]$. Now we set
    \[
        \Phi_A(x(t),y(t)) = (q(t), \ov{x}(t),\ov{y}(t),\xi(t)).
    \]
    First we note that this assignment is functorial in $A$ (see \cref{ss:weierstrass}). Clearly
    \[
        \Phi_A(x(t),y(t)) \equiv (t^d, \ov{x}_0(t),\ov{y}_0(t), \xi_0(t)) \mod \fm.
    \]
    Moreover, we have $\delta(x(t),\ov{y}(t)) \in q A[[t]]$ and thus, by \cref{l:weierstrass-polynomial}, we get $\delta(\ov{x}(t),\ov{y}(t)) \in q A[t]$. Finally,
    \begin{multline*}
        0 = f(x(t),y(t))  \equiv f(\ov{x}(t),y(t)) \\
        \equiv f(\ov{x}(t),\ov{y}(t)) + Df(\ov{x}(t),\ov{y}(t)) \cdot q(t) \theta(t) \mod q^2 A[[t]]^m.
    \end{multline*}
    Multiplying with $\Ad(Df(\ov{x}(t),\ov{y}(t)))$ yields
    \begin{multline*}
        0 \equiv \Ad(Df(\ov{x}(t),\ov{y}(t))) \cdot f(\ov{x}(t),\ov{y}(t)) + \delta(\ov{x}(t),\ov{y}(t)) q(t) \theta(t) \\
        \equiv \Ad(Df(\ov{x}(t),\ov{y}(t))) \cdot f(\ov{x}(t),\ov{y}(t)) \mod q^2 A[[t]]^m.
    \end{multline*}
    Using \cref{l:weierstrass-polynomial} again we get
    \[
        \Ad(Df(\ov{x}(t),\ov{y}(t))) \cdot f(\ov{x}(t),\ov{y}(t)) \in q^2 A[t]^m.
    \]
    Thus $\Phi_A$ is well-defined. It remains to show that $\Phi_A$ is bijective. Let $\widetilde{\gamma}$ be an $A$-valued deformation of $(Z\times\cP_\infty^n,\gamma_K)$ given by $(q(t),\ov{x}(t),\ov{y}(t),\xi(t))$. Set $x(t) \coloneqq \ov{x}(t) + q^2(t)\xi(t)$, then it suffices to show that there exists unique $\theta(t) \in A[[t]]^m$ such that $y(t) \coloneqq \ov{y}(t) + q(t)\theta(t)$ satisfies $y(t) \equiv y_0(t) \mod \fm [[t]]^m$ and $f(x(t),y(t)) = 0$. We introduce variables $Y = (Y_1,\ldots,Y_m)$ and consider the system of equations
    \begin{multline}
        \label{eq:drinfeld1}
        0 = \Ad(Df(x(t),\ov{y}(t))) \cdot f(x(t),\ov{y}(t) + q(t)Y) \\
        = \Ad(Df(x(t),\ov{y}(t))) \cdot f(x(t),\ov{y}(t)) + \delta(x(t),\ov{y}(t)) q(t) Y + q(t) \widetilde{Q}(t,Y)
    \end{multline}
    where $\widetilde{Q}(t,Y) \in A[[t]][Y]^m$ is at least quadratic in the variables $Y$. From the definition of the scheme $Z$ in \cref{l:formal-model-scheme} we have
    \begin{gather*}
        \delta(x(t),\ov{y}(t)) = q(t) v(t), \\
        \Ad(Df(x(t),\ov{y}(t))) \cdot f(x(t),\ov{y}(t)) \in q^2 A[[t]],
    \end{gather*}
    where $v(t) \in A[[t]]$ is a unit. Dividing \eqref{eq:drinfeld1} by $q^2(t)$ and multiplying with $v(t)^{-1}$ we obtain an equation of the form
    \begin{equation}
        \label{eq:drinfeld2}
        0 = P(t) + Y + Q(t,Y)
    \end{equation}
    with $P(t) \in A[[t]]^m$ and $Q(t,Y) \in A[[t]][Y]^m$ at least quadratic in $Y$. We can lift $\theta_0(t) \in K[[t]]^m$ to elements in $A[[t]]^m$ and by abuse of notation denote these by $\theta_0(t)$ as well. Then the system of equations \eqref{eq:drinfeld2} has a solution $Y = \theta_0(t) \in A[[t]]^m$ modulo $\fm$. By induction, assume we have a solution $\theta_i(t) \in A[[t]]^m$ modulo $\fm^i$. Write
    \[
        \varepsilon(t) \coloneqq - (P(t) + \theta_i(t) + Q(t,\theta_i(t))) \in \fm^i[[t]]^m.
    \]
    Then $\theta_{i+1}(t) \coloneqq \theta_i(t) + \varepsilon(t) \in A[[t]]^m$ is a solution of \eqref{eq:drinfeld2} modulo $\fm^{i+1}$. This finishes the proof that $\Phi_A$ is bijective.
\end{proof}

\begin{remark}
    Instead of the last argument, one could alternatively use the multivariate Hensel lemma for $A[[t]]$ to show that the simple root $Y = \theta_0(t)$ of $f(x(t),\ov{y}(t) + q(t)Y)$ modulo $\fm$ can be lifted to a exact solution $Y = \theta(t)$.
\end{remark}

\subsection{Formal neighborhoods via deformations}
\label{ss:formal-neighborhood}

To recap, our goal is to prove the existence of an isomorphism between the formal neighborhood of $X_\infty$ at $\alpha$ and the formal neighborhood of $Z \times \cP_\infty^n$ at $\mu(\alpha)$ up to appropriate change of coefficient field. If $\alpha$ is a $k$-rational point of $X_\infty$ then this is the original statement in \cite{Dri02}, see \cref{r:dgk-k-rational}. In the general case one has to take into account the residue field extension induced by $\mu$; the crucial point here is to use results from \cite{CdFD24} for general linear projections.

Recall that we assumed $X = V(f_1,\ldots,f_m) \subset \A^{n+m}$ where $f_1,\ldots,f_r \in k[x,y]$ with $x = (x_1,\ldots,x_n)$, $y = (y_1,\ldots,y_m)$. Moreover $\alpha \in X_\infty^{\delta,d}$ with $\delta = \det\left(\frac{\partial f_i}{\partial y_j}\right)_{i,j\leq m}$. Set $Y \coloneqq \A^n$ and consider the morphism
\[
    f \colon X \to Y,\: (x,y) \mapsto x
\]
which is the restriction of a linear projection. Write $p = \alpha(\eta)$, then by assumption we have that there exists an open neighborhood $U \subset X$ of $p$ such that the restriction of $f$ to $U$ is finite unramified. Writing $q = f(p)$, we get that the residue field extension $k_p/k_q$ is finite separable. By \cite[Theorem 3.1(2)]{CdFD24} we get that $k_\alpha/k_\beta$ is finite separable, where $\beta = f_\infty(\alpha)$ and $f_\infty \colon X_\infty \to Y_\infty$ the induced morphism on arc spaces. We use this to compare the residue fields of $\alpha$ and $\mu(\alpha)$.

Recall that $Z$ is a closed subscheme of $W = \cQ_{<d} \times \cP_{<2d}^n \times \cP_{<d}^m$. We define the morphism $\lambda \colon Z \times \cP^\infty \to Y_\infty$ as the restriction of
\[
    W \times \cP_\infty^n \to Y_\infty, \: (q,\ov{x},\ov{y};\xi) \mapsto \ov{x} + q^2 \xi. 
\]
Then it follows immediately that $\mu$ and $\lambda$ fit into a commutative diagram
\[
    \begin{tikzcd}
        X_\infty^{\delta,d} \ar[r,hookrightarrow] \ar[d,"\mu"] & X_\infty \ar[d,"f_\infty"] \\
        Z \times \cP_\infty^n \ar[r, "\lambda"] & Y_\infty.
    \end{tikzcd}
\]
In particular, $\gamma \coloneqq \mu(\alpha)$ and $\alpha$ both map to $\beta \in Y_\infty$. Now, since $k_\alpha/k_\beta$ is finite separable, so is the intermediate extension $k_\alpha/k_\gamma$. Putting together all the pieces we get the promised extension of the main result in \cite{Dri02}, which in turn finishes the proof of \cref{intro:dgk}.

\begin{theorem}
    \label{t:drinfeld}
    Let $\mu \colon X_\infty^{\delta,d} \to Z \times \cP_\infty^n$ and $\alpha \in X_\infty^{\delta,d}$ with residue field $k_\alpha$. Write $\gamma \coloneqq \mu(\alpha) \in Z \times \cP_\infty^n$ with residue field $k_\gamma$, and choose a coefficient field $\iota \colon k_\gamma\to\widehat{\cO_{Z\times \cP_\infty^n,\gamma}}$. Then there exists an isomorphism of formal neighborhoods
    \[
        (X_\infty, \alpha) \isom (Z \times \cP_\infty^n, \gamma) \times_{k_\gamma} k_\alpha,
    \]
    where on the right hand side we consider the fiber product in the category of formal schemes.
\end{theorem}

\begin{proof}
    We write $S \coloneqq \cO_{X_\infty, \alpha}$ with maximal ideal $\fn$ and residue field $k_\alpha$. Similarly, write $R \coloneqq \cO_{Z \times \cP_\infty^n,\gamma}$ with maximal ideal $\fm$ and residue field $k_\gamma$. The morphism $\mu$ determines a field extension $j \colon k_\gamma \to k_\alpha$ and a $k_\alpha$-point $\gamma_{k_\alpha}$ of $Z\times \cP_\infty^n$. By \cref{t:deformations} we have a natural isomorphism
    \[
        \Def_{S} \isom \Def_{(X_\infty,\alpha)} \isom \Def_{(Z \times \cP_\infty^n,\gamma_{k_\alpha})} \isom \Def_{R,k_\alpha/k_\gamma}.
    \]
    Write $\cR = (R/\fm^i)_i,\: \cS = (S/\fn^i)_i \in \Cpt_K$. The choice of coefficient field $\iota$ for $\^R$ induces a coefficient field for $\cR$ and thus by \cref{p:deformations-finite-separable} an isomorphism of pro-objects $\cR \otimes_{k_\gamma} k_\alpha \isom \cS$. By \cref{l:base-change-cpt} we get an isomorphism in $\TopRing_k$
    \[
        \^{\cO_{Z\times\cP_\infty^n,\gamma}} \^\otimes_{k_\gamma} k_\alpha \isom \^{\cO_{X_\infty,\alpha}}.
    \]
    The statement now follows by applying the functor $\Spf$.
\end{proof}

\begin{remark}
    \label{r:dgk-k-rational}
    Note that if $\alpha \in X_\infty$ is $k$-rational, then so are $\beta$ and $\gamma$. In particular, there is only one choice of coefficient field in \cref{t:drinfeld} and the resulting isomorphism is $(X_\infty,\alpha) \isom (Z \times \cP_\infty^n,\gamma)$, thus fully recovering the original statement in \cite{Dri02}.
\end{remark}

\section{Application to embedding codimension}
\label{s:embcodim}

\subsection{On embedding codimension}
\label{ss:on-embcodim}

We recall the definition of embedding codimension in \cite{CdFD22}. If $(A,\fm)$ is a local ring, then the embedding codimension of $A$ was defined there as
\[
    \embcodim(A) \coloneqq \height(\ker \gamma),
\]
where $\gamma \colon \Sym_K(\fm/\fm^2) \to \gr(A)$ is the natural surjection and $K$ denotes the residue field of $A$. If $A$ is Noetherian, then $\dim A = \dim \gr(A)$ and hence
\[
    \embcodim(A) = \embdim(A) - \dim (A).
\]
Our aim is ultimately to compare embedding codimensions using \cref{intro:dgk}. The statement there involves a base change of the completion of a local ring. In the non-Noetherian case it is not clear whether this base change is again the completion of a local ring, as for non-Noetherian completions the inverse limit topology does not coincide with the preadic topology as a local ring, see \SPC{05JA}. Hence we need to extend the notion of embedding codimension in that case. We first define a sufficiently large subcategory of $\Cpt_K$ resp.\ of $\TopRing_k$.

\begin{definition}
    Let $\cR = (R_i)_i \in \Cpt_K$ and $\^R = \varprojlim_i R_i$ the limit in $\TopRing_k$ with maximal ideal $\fm$. Then $\^R$ (or $\cR$) is called \emph{quasi-adic} if the closures $\^{\fm^n}$ of the ideals $\fm^n$ form a basis for the topology of $\^R$.
\end{definition}

\begin{remark}
    \label{r:quasi-adic-surjection}
    If $\cR$ is quasi-adic, then we consider $\^R$ with the filtration given by the ideals $\^{\fm^n}$. Then the natural map
    \[
        \^\gamma \colon \Sym_K(\fm/\^{\fm^2}) \to \gr(\^R)
    \]
    is surjective. Indeed, this follows from the fact that the maps
    \[
        \fm^n / \^{\fm^{n+1}} \to \^{\fm^n}/\^{\fm^{n+1}}
    \]
    are bijective. Conversely, this property characterizes quasi-adic rings, see \cite{Chi21}.
\end{remark}

\begin{remark}
    If $(R,\fm)$ is any local ring, then the completion $\^R$ is quasi-adic. Moreover, if $\^R$ is quasi-adic and the maximal ideal $\^\fm$ is finitely generated, then $\^R$ is adic. This follows essentially from \SPC{09B8} (see also \cite{Chi21}). Conversely, for $R$ the localization of $k[x_i \mid i\in I]$ at the maximal ideal $\fm = (x_i \mid i \in I)$, the completion of $R$ is adic if and only if $I$ is finite \SPC{05JA}.
\end{remark}

\begin{example}
    For an example of $\cR \in \Cpt_K$ such that $R = \varprojlim_i R_i$ is not quasi-adic, let $(A,\fm)$ be a quasi-adic but not adic ring with residue field $K$. Choose a coefficient field $\iota \colon K \to A$. Consider the ring of restricted power series in $A$
    \[
        A\{t\} \coloneqq \{f(t) \in A[[t]] \mid f(t) = \sum_{i\geq0} f_i t^i, f_i \in \^{\fm^i}\},
    \]
    where we set $\^{\fm^0} \coloneqq \iota(K)$. We regard $A\{t\}$ as a linear topological ring with basis for the topology given by
    \[
        \fn_i \coloneqq \{f(t) \in A\{t\} \mid \ord_t f(t) \geq i\}.
    \]
    It is straightforward to check that $A\{t\}$ is complete with respect to this topology. Equivalently, the assignment $\cR(i) \coloneqq A\{t\}/\fn_i$ defines an object of $\Cpt_K$. We claim that $A\{t\}$ is not quasi-adic. Indeed, as $A$ is not adic we can find $n\geq 0$ and for each $i\geq n$ an element $a_i \in \^{\fm^i}$ with $a_i \notin \fm^n$. Then $a_i t^i \in \fn_i$, but $a_i t^i \notin \fn_1^n + \fn_{n+1}$. This shows that the closure of $\fn_1^n$ is not open.
\end{example}

\begin{definition}
    Let $\cR \in \Cpt_K$ be quasi-adic. We define the \emph{embedding codimension} of $\cR$ as
    \[
        \embcodim(\cR) \coloneqq \height(\ker \^\gamma ),
    \]
    where $\^\gamma \colon \Sym_K(\fm/\^{\fm^2}) \to \gr(\^R)$ is the map from \cref{r:quasi-adic-surjection}.
\end{definition}

For local rings this definition is indeed consistent in the following sense.

\begin{lemma}
    If $(R,\fm)$ is a local ring with residue field $K$ and $\cR$ denotes the corresponding object in $\Cpt_K$ defined by $\cR(i) \coloneqq R/\fm^{i+1}$, then we have
    \[
        \embcodim(R) = \embcodim(\cR).
    \]
\end{lemma}

\begin{proof}
    Write $S$ for the $\fm$-adic completion of $R$ and $\fn \coloneqq S \to R/\fm$ for the maximal ideal of $S$. For all $i\geq 1$ we have
    \[
        \fm^i/\fm^{i+1} \isom \^{\fn^i}/\^{\fn^{i+1}}
    \]
    and the statement follows.
\end{proof}

We make the following observations.

\begin{lemma}
    \label{l:embcodim-base-change}
    Let $\cR \in \Cpt_K$ be quasi-adic. Let $K \to \cR$ be a coefficient field. If $K \subset L$ is a field extension, then $\cR \otimes_K L$ is again quasi-adic and
    \[
        \embcodim(\cR) = \embcodim(\cR \otimes_K L).
    \]
\end{lemma}

\begin{proof}
    Write $R \coloneqq \varprojlim_i R_i$ with maximal ideal $\fm$. Then $R_L \coloneqq \varprojlim_i (R_i \otimes_K L)$ is again quasi-adic with maximal ideal $\fm_L$. If $\^\gamma$ is the natural surjection as in \cref{r:quasi-adic-surjection}, then the corresponding map $\^\gamma_L$ is just the base change of $\^\gamma$ to $L$.
\end{proof}

Recall that a ring homomorphism $R \to S$ is called essentially smooth if it is the localization of a smooth ring homomorphism or equivalently, if it is formally smooth and essentially finitely presented.

\begin{proposition}
    \label{p:embcodim-formally-smooth-extension}
    Let $(R,\fm)$ a local ring essentially of finite type over $k$. Let $(R,\fm) \to (S,\fn)$ be a local $k$-algebra map that is the direct limit of essentially smooth local $k$-algebra maps $(R,\fm) \to (S_i,\fn_i)$, with essentially smooth transition maps $(S_i,\fn_i) \to (S_j,\fn_j)$. Then
    \[
        \embcodim(S) = \embcodim(R).
    \]
\end{proposition}

\begin{proof}
    Write $K = R/\fm$, $L = S/\fn$ and $L_i = S_i/\fn_i$. We use an argument similar to \cite[Theorem 8.3]{CdFD22} to show that
    \[
        \embcodim(S) = \limsup_{i\in \N} \embcodim(S_i).
    \]
    Let us sketch the proof. Since $(S_i,\fn_i) \to (S_j,\fn_j)$ is essentially smooth, we have a commutative diagram
    \[
        \begin{tikzcd}
            \Sym_{L_j}(\fn_j/\fn_j^2) \otimes_{L_j} L \ar[r, "\gamma_j"] & \gr(S_j) \otimes_{L_j} L \\
            \Sym_{L_i}(\fn_i/\fn_i^2) \otimes_{L_i} L \ar[r, "\gamma_i"] \ar[u, hookrightarrow] & \gr(S_i) \otimes_{L_i} L \ar[u],
        \end{tikzcd}
    \]
    where the left vertical arrow is an extension of polynomial rings and hence faithfully flat. The same holds when replacing $\gamma_j$ with the surjection $\gamma \colon \Sym_L(\fn/\fn^2) \to \gr(S)$. Setting $\fa \coloneqq \ker \gamma$ and $\fa_i \coloneqq \ker \gamma_i$ one has $\fa = \varinjlim_i \fa_i$ and hence $\height (\fa) = \limsup_i \height (\fa_i)$.
    
    Now, for all $R \to S_i$ we have
    \[
        \embdim(S_i) = \embdim(R) + \embdim(S_i \otimes_R K)
    \]
    and
    \[
        \dim(S_i) = \dim(R) + \dim(S_i \otimes_R K).
    \]
    Since $R \to S_i$ is essentially smooth the fiber $S_i \otimes_R K$ is regular. Thus $\embcodim(R) = \embcodim(S)$.
\end{proof}

It is important to state here that even for schemes of finite type over $k$, the embedding codimension does not satisfy any obvious semicontinuity properties, as the following example shows.

\begin{example}
    \label{ex:embcodim-not-semicontinuous}
    Consider the scheme $X \subset \A_\C^3$ defined by the product of the ideals $(x)$ and $(y,z)^2$. Geometrically $X$ is the union of a plane $X_1$ and a double line $X_2$. The locus of points of embedding codimension $\leq 1$ is the closed subset $X_1$, which consists of the (open) smooth locus $X_1 \setminus 0$ and the origin as the only point of embedding codimension $1$. In particular, the generic point $\eta$ of $X_2$ is a generization of $0$ such that
    \[
        \embcodim(\cO_{X,0}) = 1 < 2 = \embcodim(\cO_{X,\eta}).
    \]
\end{example}

In what follows we prove that the drop in the minimal dimension of a component observed in \cref{ex:embcodim-not-semicontinuous} is indeed the only obstruction to the embedding codimension being upper semi-continuous. We assume the statement is well-known to experts, but were not able to find a reference in the literature.

\begin{proposition}
    \label{p:emcodim-generically-semicontinuous}
    Let $X$ be a scheme of finite type over $k$. Let $X' \subset X$ be a closed irreducible subset. Then there exists a nonempty open subset $U$ of $X'$ such that the function
    \[
        x \mapsto \embcodim(\cO_{X,x})
    \]
    is upper semi-continuous when restricted to $U$. In particular, there exists a nonempty open subset $U'$ of $X'$ where $\embcodim(\cO_{X,x})$ is constant.
\end{proposition}

\begin{proof}
    We may assume that $X = \Spec S$, where $S = R/I$ with $R = k[x_1,\ldots,x_n]$. Let $\eta$ be the generic point of $X'$; it corresponds to a prime $\fp$ of $S$ resp.\ of $R$. By \cite[0, (14.2.6)]{egaiv1}, the function $\fq \mapsto \height (I R_\fq)$ is lower semi-continuous on $X'$. Hence there exists a nonempty open subset $U$ of $X'$ such that for all $\fq \in U$ we have $\height(I R_\fq) = \height(I R_\fp) = r$.

    Now let $\fq \in \Spec S$ with residue field $L$. Since $k$ is perfect, the conormal sequence
    \[
        \begin{tikzcd}[column sep = small]
            0 \ar[r] & \fq S_\fq/\fq^2 S_\fq \ar[r] & \Om_{S/k} \otimes L \ar[r] & \Om_{L/k} \ar[r] & 0
        \end{tikzcd}
    \]
    is short exact. Then $\dim_L \Om_{L/k} = \trdeg_k L = \dim S/\fq$ and we have
    \[
        \dim_L \Om_{S/k} \otimes L = \embdim(S_\fq) + \dim S/\fq. 
    \]
    Moreover, we have
    \[
        \dim S/\fq = n - \height(I R_\fq) - \dim S_\fq.
    \]
    Now let $I = (f_1,\ldots,f_s)$ and $Df \coloneqq (\frac{\partial f_i}{\partial x_j})_{i,j}$. We write $Df(\fq)$ for $Df$ evaluated at $\fq$. Take the conormal sequence for the surjection $R \to S$ and basechange to $L$ to get
    \[
        \begin{tikzcd}[column sep = small]
            I/I^2 \otimes L \ar[r] & \Omega_{R/k} \otimes L \ar[r] & \Omega_{S/k} \otimes L \ar[r] & 0.
        \end{tikzcd}
    \]
    Since $\Om_{S/k} \otimes L$ is the cokernel of $Df(\fq)$, we get
    \[
        \dim_L (\Om_{S/k} \otimes L) = n - \rk(Df(\fq)).
    \]
    Putting everything together we have
    \[
        \embcodim(S_\fq) = \height(I R_\fq) - \rk(Df(\fq)),
    \]
    which on $U$ is the difference between a constant function and a lower semi-continuous one.
\end{proof}

\subsection{Embedding codimension over maximal divisorial sets}
\label{ss:max-div}

Let us first give the proof of \cref{intro:embcodim}. As in \cref{ss:complete-intersection}, we reduce first to $X$ affine and then, using \cref{p:complete-intersection}, to the case where $X$ is a complete intersection. That way it suffices to prove the following proposition.

\begin{proposition}
    Let $X = V(f_1,\ldots,f_m) \subset \A^{n+m}$ as in \cref{s:drinfeld}, with $f_i \in k[x,y]$. Let
    \[
        \delta \coloneqq \det \left(\frac{\partial f_i}{\partial y_j}\right)_{i,j\leq m}
    \]
    and let $\beta \in X_\infty$ with $\ord_\beta \delta < \infty$. Write $W \coloneqq \overline{\{\beta\}} \subset X_\infty$. Then there exists an open subset $W^0 \subset W$ such that the function
    \begin{equation}
        \label{eq:embcodim-fct-1}
        W \to \N, \: \alpha \mapsto \embcodim(\cO_{X_\infty,\alpha})
    \end{equation}
    is constant on $W^0$.
\end{proposition}

\begin{proof}
    Let $d = \ord_\beta \delta$ and recall that $X_\infty^{\delta,d}$ is the locally closed subset defined by
    \[
        X_\infty^{\delta,d} \coloneqq \{\alpha \in X_\infty \mid \ord_\alpha \delta = d\}.
    \]

    Let $W^{\delta,d} \coloneqq W \cap X^{\delta,d}_\infty$. Note that $W^{\delta,d}$ is an open subset of $W$. Consider the morphism $\mu \colon X_\infty^{\delta,d} \to Z \times \cP_\infty^n$ from \cref{ss:formal-model}. Write
    \begin{equation}
        \label{eq:proj-formal-model}
        \mu_Z \colon X_\infty^{\delta,d} \to Z \times \cP_\infty^n \to Z
    \end{equation}
    for the composition of $\mu$ with the projection to $Z$.  We define the following function:
    \begin{equation}
        \label{eq:embcodim-fct-2}
        W^{\delta,d} \to \N,\: \alpha \mapsto \embcodim(\cO_{Z,\mu_Z(\alpha)}).
    \end{equation}
    We first claim that the function \eqref{eq:embcodim-fct-2} equals the restriction of \eqref{eq:embcodim-fct-1} to $W^{\delta,d}$. Indeed, by \cref{t:drinfeld,l:embcodim-base-change} we have
    \[
        \embcodim(\cO_{X_\infty,\alpha}) = \embcodim(\cO_{Z\times \cP_\infty^n,\mu(\alpha)})
    \]
    for each $\alpha \in X_\infty^{\delta,d}$. By \cref{p:embcodim-formally-smooth-extension} it follows that
    \[
        \embcodim(\cO_{Z\times \cP_\infty^n,\mu(\alpha)}) = \embcodim(\cO_{Z,\mu_Z(\alpha)}).
    \]
        
    Let $Z_W$ denote the closure of $\mu_Z(\alpha)$ inside $Z$. By \cref{p:emcodim-generically-semicontinuous} there exists a nonempty open subset $U_W$ of $Z_W$ such that the function $z \mapsto \embcodim(\cO_{Z,z})$ is constant on $U_W$. Define $W^0 \coloneqq W \cap \mu_Z^{-1}(U_W)$. Then the restriction of \eqref{eq:embcodim-fct-1} to $W^0$ is constant.
\end{proof}

\begin{remark}
    \label{r:embdim-formal-model}
    \cref{intro:embcodim} can be seen as an extension of \cite[Theorem 10.5]{CdFD24}. To briefly summarize the argument given there, let $\mu \colon X_\infty^{\delta,d} \to Z \times \cP_\infty$ be the morphism from \cref{ss:formal-model}. It is shown that for any $k$-rational $\alpha \in X_\infty^{\delta,d}$ we have that $\embdim( \cO_{Z,\mu_Z(\alpha)})$ is constant. Therefore the function
    \[
        \alpha \in X_\infty^{\delta,d}(k) \to \embcodim(\cO_{X_\infty,\alpha})
    \]
    is the difference of a constant and an upper-semicontinuous function, and hence lower semi-continuous itself. This also suggests that, to control the embedding codimension of $\alpha \in X_\infty^{\delta,d}$, it may suffice to control the local dimension of $Z$ at the image of $\mu$.
\end{remark}

We now want to detail how \cref{intro:embcodim} relates to invariants of divisorial valuations on a variety $X$. Recall that a valuation $\nu$ of the function field $k(X)$ with values in $\Z$ and center in $X$ is called divisorial if its residue field $k_\nu$ has transcendence degree $\dim X - 1$ over $k$. Equivalently, $\nu$ is of the form $\nu = q \ord_E$ where $q \in \Z_{>0}$ and $E$ is a prime divisor on $Y$ normal with $f \colon Y \to X$ proper birational. In this way one defines the following variants of the discrepancy of $\nu$.

\begin{definition}
    For $\nu = q \ord_E$ a divisorial valuation as above, we define
    \begin{enumerate}
        \item the \emph{Mather log discrepancy} to be 
        \[
            \^a_\nu(X) \coloneqq q (\ord_E(\Jac_f) + 1),
        \]
        \item and the \emph{Mather-Jacobian log discrepancy} to be
        \[
            a^{\MJ}_\nu(X) \coloneqq q(\ord_E(\Jac_f) - \ord_E(\Jac_X) + 1).
        \]
    \end{enumerate}
\end{definition}

Mather discrepancies featured prominently in the change-of-variables formula in motivic integration \cite{DL99} and were further studied in \cite{dFEI08,Ish13}, whereas Mather-Jacobian discrepancies first appeared in \cite{dFD14,EIM16}. If $X$ is in addition $\Q$-Gorenstein and $a_\nu(X)$ denotes the usual discrepancy of $\mu$, one has the relations
\[
    a^{\MJ}_\nu(X) \leq a_\nu(X) \leq \^a_\nu(X),
\]
with the first being an equality when $X$ is a local complete intersection, and the second when $X$ is smooth \cite[Section 3.2]{dFD14}.

Mather(-Jacobian) discrepancies are intrinsically linked to the arc space, with this relation usually formulated in terms of cylindrical subsets associated to each divisorial valuation as follows.

\begin{definition}
    Let $X$ be a variety over $k$ and $\nu$ a divisorial valuation on $X$. The maximal divisorial set associated to $\nu$ is the subset of $X_\infty$ defined by
    \[
        C_\nu(X) \coloneqq \ov{\{\alpha \in X_\infty \mid \ord_\alpha = \nu\}}.
    \]
\end{definition}

As before, let $f \colon Y \to X$ be proper birational with $Y$ normal and such that $\nu = q \ord_E$ for a prime divisor $E$ on $Y$. We write $E^0 \subset E$ to be the open subset of $E$ where $E, Y$ are smooth and no other component of the exceptional locus intersects $E$. Then by \cite[Lemma 11.3]{dFD20} we have
\[
    C_\nu(X) = \ov{\Cont^{\geq q}(E^0,Y)}.
\]
In particular $C_\nu(X)$ is irreducible with generic point $\alpha_\nu$; we call $\alpha_\nu$ the \emph{maximal divisorial arc} associated to $\nu$. In fact, $C_\nu(X)$ is what is often called a \emph{cylindrical subset}; that is, it is of the form $\pi_n^{-1}(V)$, where $\pi_n \colon X_\infty \to X_n$ and $Z \subseteq X_n$ is constructible. For cylindrical subsets one can define a notion of codimension, and this codimension of $C_\nu(X)$ inside $X_\infty$ equals $\^a_\nu(X)$ \cite[Theorem 3.8]{dFEI08}. Alternatively, one has the the following result, relating Mather(-Jacobian) log discrepancies to invariants of the local ring at $\alpha_\nu$.

\begin{theorem}[{\cite[Theorem 11.5]{CdFD24}}]
    \label{t:discrepancies-max-div-arc}
    Let $X$ be a variety over a perfect field $k$, $\nu$ a divisorial valuation on $X$ and $\alpha_\nu$ the corresponding maximal divisorial arc. Then
    \begin{enumerate}
        \item \label{eq:discrepancies-1} $\embdim (\cO_{X_\infty,\alpha_\nu}) = \^a_\nu(X)$, and
        \item \label{eq:discrepancies-2} $\dim (\^{\cO_{X_\infty,\alpha_\nu}}) \geq a_\nu^{\MJ}(X)$.
    \end{enumerate}
\end{theorem}

In the case where $k$ is of characteristic $0$, \cref{t:discrepancies-max-div-arc} was first proven in \cite{MR18}. Let us sketch the proof in the general case. The equality in \eqref{eq:discrepancies-1} is deduced by using a version of the birational transformation rule, expressed in terms of the embedding dimension of maximal divisorial arcs \cite[Theorem 9.2]{dFD20}. The inequality in \eqref{eq:discrepancies-2} is then an immediate consequence of \eqref{eq:discrepancies-1} and the following general bound on the embedding codimension.

\begin{theorem}[{\cite[Theorem 9.8]{CdFD24}}]
    \label{t:embcodim-bound}
    Let $X$ be a scheme locally of finite type over a perfect field $k$. For any $\alpha \in X_\infty$, we have $\alpha(\eta) \in X$ is smooth if and only if
    \[
        \embcodim(\cO_{X_\infty, \alpha}) \leq \ord_\alpha(\Jac_{X^0}) < \infty,
    \]
    where $X^0$ is the unique irreducible component of $X$ containing $\alpha(\eta)$.
\end{theorem}

In contrast, \cref{intro:embcodim} applied to the maximal divisorial subset $C_\nu(X)$ together with \cref{p:embcodim-formally-smooth-extension} gives the following.

\begin{corollary}
    \label{c:embcodim-max-div}
    There exists an nonempty open subset $C_\nu(X)^0$ of $C_\nu(X)$ such that the function
    \[
        \alpha \mapsto \embcodim(\cO_{X_\infty, \alpha})
    \]
    is finite constant on $C_\nu(X)$. In particular, for any $\alpha \in C_\nu(X)^0$ the embedding codimension of a finite formal model for $\alpha$ equals $\embcodim (\cO_{X_\infty,\alpha_\nu})$.
\end{corollary}

We emphasize that the explicit bound in \cref{t:embcodim-bound} does not immediately follow from \cref{c:embcodim-max-div}, as discussed in \cite[Section 10]{CdFD22}. 

We anticipate that the strategy of considering the scheme of formal models as in \cref{t:drinfeld} will yield further results on invariants of singularities of the arc space. However, we want to emphasize that this is not straightforward even when trying to find a similar relation for the (embedding) dimension. As noted in the introduction, the first obstacle is that both dimension and embedding dimension obviously depend on the choice of finite formal model for $\alpha \in C_\nu(X)$. One may circumvent this by considering the \emph{minimal formal model} instead: for any $\alpha \in X_\infty(k) \setminus (\Sing X)_\infty$ there exists a scheme $Z$ of finite type such that
\[
    (X_\infty,\alpha) \isom (Z,z) \times (\A^\N,0)
\]
and $(Z,z)$ itself is not of the form $(Z,z) \isom (Z',z') \times (\A^1,0)$. Note that the formal scheme $(Z,z)$ is unique up to isomorphism. Unfortunately, as was observed in \cite[Section 6]{BMC24}, for a divisorial valuation $\nu$ on a curve $X$ both dimension and embedding dimension of the minimal formal model of a general $k$-rational arc in $C_\nu(X)$ are strictly smaller than those of $\cO_{X_\infty,\alpha_\nu}$. As the completion of $\cO_{X_\infty,\alpha_\nu}$ is Noetherian however, we can similarly write
\[
    (X_\infty,\alpha_\nu) \isom \^Z_\nu \times (\A^m_{k_{\alpha_\nu}},0)
\]
with $\^Z_\nu$ a Noetherian formal scheme unique up to isomorphism and not of the form $\^Z_\nu \isom \^Z' \times (\A^1_{k_{\alpha_\nu}},0)$. We call $\^Z_\nu$ the \emph{minimal formal model} of $(X_\infty,\alpha_\nu)$. A follow-up question worth investigating is thus:

\begin{question}
    Let $X$ be a variety over a perfect field $k$ and $\nu$ a divisorial valuation on $X$. Denote by $\alpha_\nu$ the maximal divisorial arc and let $\alpha$ be a general $k$-rational arc in $C_\nu(X)$. If $(Z,z)$ and $\^Z_\nu$ are the minimal formal models of $(X_\infty,\alpha)$ and $(X_\infty, \alpha_\nu)$ respectively, do we have $\embdim(\cO_{Z,z}) = \embdim(\cO_{\^Z_\nu})$ (and similar for the dimensions)?
\end{question}

Following the observations in \cref{r:embdim-formal-model}, we hope that a closer study of the geometry of the scheme of formal models will eventually provide an answer to the above question, as well as more insight on the singularities of maximal divisorial sets more generally.

\subsection{On formally smooth arcs}

If $(R,\fm)$ is a Noetherian local ring, then $\dim \gr(R) = \dim R$ and thus $\embcodim(R) = 0$ if and only if $R$ is regular (which, since $k$ is perfect, is equivalent to $R$ being essentially smooth over $k$). In the general case we have the following result.

\begin{theorem}[{\cite[0, (19.5.4)]{egaiv1}}]
    \label{t:formal-smoothness-embcodim}
    Let $(R,\fm)$ be a local ring over a perfect field $k$. Then $\embcodim(R) = 0$ if and only if $R$ is formally smooth over $k$ for the $\fm$-adic topology. That is, for any $k$-algebra $C$, ideal $J \subset C$ with $J^2 = 0$ and diagram
    \[
        \begin{tikzcd}
            R \ar[r, "\ov{\tau}"] \ar[rd, dashed, "\tau"] & C/J \\
            k \ar[r] \ar[u] & C \ar[u]
        \end{tikzcd}
    \]
    with $\ker(\ov{\tau}) \supset \fm^n$ for some $m\geq 0$ there exists a lift $\tau \colon R \to C$ making the diagram commute.
\end{theorem}

\cref{t:formal-smoothness-embcodim} makes precise the idea that the embedding codimension is an invariant of singularities in infinite dimensions. This begs the question whether there exists a characterization of those arcs whose local ring has embedding codimension $0$. In the case of $k$-rational arcs this question was resolved in \cite{BS17b}, whose main result is the following.

\begin{theorem}[{\cite[Theorem 1.6]{BS17b}}]
    \label{t:bs-smooth-rational-arcs}
    Let $X$ be a variety over a perfect field $k$. Let $\alpha \in X_\infty(k)$ be a $k$-rational arc. Then the following are equivalent:
    \begin{enumerate}
        \item \label{eq:bs-smooth-rational-arcs} There exists unique smooth formal branch at $\alpha(0)$ containing the image of $\alpha$, and
        \item The local ring $\cO_{X_\infty,\alpha}$ is formally smooth for the $\fm_\alpha$-adic topology, or equivalently,
        \[
            \embcodim(\cO_{X_\infty}) = 0.
        \]
    \end{enumerate}
\end{theorem}

To elaborate on statement \eqref{t:bs-smooth-rational-arcs}, it was proven in \cite[Proposition 3.6]{BS17b} that for any $k$-rational arc $\alpha$ which is not contained in $\Sing X$ there exists a unique minimal prime $\fq$ of $\widehat{\cO_{X,\alpha(0)}}$ such that $\alpha$ factors through $\widehat{\cO_{X,\alpha(0)}}/\fq$. The assertion that the unique formal branch containing $\alpha$ is smooth just means that $\widehat{\cO_{X,\alpha(0)}}/\fq$ is formally smooth over $k$. In the case $
\alpha \in (\Sing X)_\infty$, then \cite[Theorem 8.7]{CdFD22} says that $\embcodim(\cO_{X_\infty,\alpha}) = \infty$. Thus the conclusion of \cref{t:bs-smooth-rational-arcs} holds in this case as well.

The main result of this section is a generalization of \cref{t:bs-smooth-rational-arcs} to arcs which are not $k$-rational. For the proof we will reduce to the $k$-rational case, and as such, we will need some mild assumption in order to ensure the existence of enough $k$-rational points. Note that even if $k$ is algebraically closed, if its cardinality is countable, then by \cite[Proposition 2.11]{Ish04} there exist closed subsets of the arc space which do not have a $k$-rational points; hence we will need to assume $k$ uncountable in general. However, for \emph{stable} points, as is the case in \cite[Question 2.10]{Reg18}, there always exist $k$-rational points for $k$ algebraically closed. Let us first recall the definition of stable points from e.g.\ \cite{Reg18}.

\begin{definition}
    \label{d:stable}
    Let $X$ be a variety over a perfect field $k$ and $\alpha \in X_\infty$. Then $\alpha$ is called \emph{stable} if there exists an $N \geq 0$, an affine open $U \subset X_N$ containing the image of $\alpha$ such that for $n \geq N$, the truncation morphisms $\pi_{n+1,n} \colon X_{n+1} \to X_n$ induce a trivial fibration
    \[
        \ov{\{\pi_{n+1}(\alpha)\}}\cap \pi_{n+1,N}^{-1}(U) \to \ov{\{\pi_{n}(\alpha)\}}\cap \pi_{n,N}^{-1}(U)
    \]
    with fiber $\A^{\dim X}$, where $\pi_n \colon X_\infty \to X_n$.
\end{definition}

By \cite[Corollary 11.5]{dFD20} any maximal divisorial arc is stable. In fact, by \cite[Theorem 10.8]{dFD20} stable points are precisely those arcs whose local rings have finite embedding dimension.

\begin{lemma}
    \label{l:stable-k-rational}
    Let $X$ be a variety over a perfect field $k$. Let $\alpha \in X_\infty$ be stable and $W$ a nonempty open of $\ov{\{\alpha\}}$. Then there exists a point $\beta \in W$ with residue field $\ov{k}$.
\end{lemma}

\begin{proof}
    We may assume there exists $N \geq 0$ and $U \subset X_N$ as in \cref{d:stable} and such that $\pi_N^{-1}(U) \subset W$. Pick any closed point $\beta_N$ of $\ov{\{\pi_N(\alpha)\}} \cap U \subset X_N$. As $\pi_{N+1,N}$ is locally a trivial fibration over $\ov{\{\pi_N(\alpha)\}} \cap U$ there exists a $\ov{k}$-rational point $\beta_{N+1} \in \ov{\{\pi_{N+1}(\alpha)\}}\cap \pi_{N+1,N}^{-1}(U) $ in the fiber over $\beta_N$. By induction this gives a $\ov{k}$-rational point in $W$.
\end{proof}

\begin{theorem}
    \label{t:smooth-arcs}
    Let $X$ be a variety over a perfect field $k$. Let $\alpha \in X_\infty$ and assume that $X$ is geometrically unibranch at $\alpha(0)$. Assume that either
    \begin{enumerate}
        \item $k$ is uncountable, or
        \item $\alpha$ is stable.
    \end{enumerate}
     Then $X$ is smooth at $\alpha(0)$ if and only if $\embcodim(\cO_{X_\infty,\alpha}) = 0$.
\end{theorem}

\begin{proof}
    Assume that $X$ is geometrically unibranch at $x \coloneqq \alpha(0)$ and $\embcodim(\cO_{X_\infty,\alpha}) = 0$. We first prove that we can reduce to the case where $k$ is algebraically closed. Denote by $\ov{k}$ the algebraic closure of $k$ and let $X_{\ov{k}}$ be the base change of $X$. Then the base change $(X_\infty)_{\ov{k}}$ is the relative arc space of $X_{\ov{k}}$ over $\ov{k}$. Pick a lift $\beta \in (X_\infty)_{\ov{k}}$ of $\alpha$. We have that $\beta(0)$ maps to $\alpha(0)$ via $X_{\ov{k}} \to X$. To conclude we need to verify the following properties:
    \begin{enumerate}
        \item \label{eq:alg-cl-1} $X_{\ov{k}}$ is geometrically unibranch at $\beta(0)$,
        \item \label{eq:alg-cl-2} $\embcodim(\cO_{(X_\infty)_{\ov{k}},\beta}) =0$,
        \item \label{eq:alg-cl-3} $\beta(0)$ smooth implies that $\alpha(0)$ is smooth, and
        \item \label{eq:alg-cl-4} if $\alpha$ is stable then so is $\beta$. 
    \end{enumerate}
    Assertion \eqref{eq:alg-cl-1} follows from \SPC{0C55}. For \eqref{eq:alg-cl-2} we first note that, by \cref{t:formal-smoothness-embcodim}, the local ring $\cO_{X_\infty,\alpha}$ is formally smooth over $k$ for the $\fm_\alpha$-adic topology. Note that $\cO_{(X_\infty)_{\ov{k}},\beta}$ is a localization of $\cO_{X_\infty,\alpha} \otimes_k \ov{k}$ and hence formally smooth over $\ov{k}$ for the $\fm_\alpha$-adic topology by \cite[0, (19.3.5)]{egaiv1}. Since the latter is finer than the $\fm_\beta$-adic topology the claim follows from \cite[0, (19.3.8)]{egaiv1}. Assertion \eqref{eq:alg-cl-3} is just a special case of the argument above. Finally, for \eqref{eq:alg-cl-4} assume that $\alpha$ is stable. By \cite[Theorem 10.8]{dFD20} this is equivalent to $\embdim(\cO_{X_\infty,\alpha}) < \infty$. Consider the local ring homomorphism
    \[
        \varphi \colon \cO_{X_\infty,\alpha} \to \cO_{(X_\infty)_{\ov{k}},\beta}.
    \]
    Since $k \subset \ov{k}$ is a direct limit of finite separable extensions it is formally \'etale and hence so is $\varphi$. Now consider the diagram
    \[
        \begin{tikzcd}[column sep = small]
            \cO_{(X_\infty)_{\ov{k}},\beta} \ar[rrd, dashed] \ar[rr, two heads] & & k_\beta \\
            \cO_{X_\infty,\alpha} \ar[u, "\varphi"] \ar[r, two heads] & \cO_{X_\infty,\alpha} / \fm_\alpha^2 \isom k_\alpha \oplus \fm_\alpha/\fm_\alpha^2 \ar[r] & k_\beta \oplus (\fm_\alpha/\fm_\alpha^2 \otimes_{k_\alpha} k_\beta) \ar[u, two heads].
        \end{tikzcd}
    \]
    By the infinitesimal lifting criterion there exists a unique diagonal arrow making the diagram commute. This in turn means that the map of $k_\beta$-vector spaces induced by $\varphi$
    \[
        \fm_\alpha/\fm_\alpha^2 \otimes_{k_\alpha} k_\beta \to \fm_\beta/\fm_\beta^2
    \]
    has a unique retraction, which implies it is an isomorphism. Therefore $\embdim(\cO_{(X_\infty)_{\ov{k}},\beta}) < \infty$ and hence $\beta$ is stable.
    
    Now assume $k$ is algebraically closed. \cref{intro:embcodim} applied to $W = \ov{\{\alpha\}}$ gives a nonempty open subset $W^0$ of $W$ where the embedding codimension is $0$. By \cite[Proposition 2.10]{Ish04} in case $k$ is uncountable, or by \cref{l:stable-k-rational} in case $\alpha$ is stable, there exists a $k$-rational arc $\alpha' \in W^0$. Then $x' \coloneqq \alpha'(0)$ is a specialization of $x = \alpha(0)$. By \cref{t:bs-smooth-rational-arcs} the unique formal branch at $x'$ containing the image of $\alpha'$ is smooth. Applying \SPC{0CB4} gives an \'etale morphism $f\colon U \to X$ and a point $u' \in U$ mapping to $x'$ such that the irreducible components of $U$ at $u'$ correspond to the branches of $X$ at $x'$. Let $\beta' \in U_\infty$ be the unique lift of $\alpha'$ with $\beta'(0) = u'$ and $U_i \subset U$ the unique component containing the image of $\beta'$. The fact that the branch containing $\alpha'$ is smooth is equivalent to $U_i$ being smooth at $u'$. Since $f$ is \'etale, so is $f_\infty \colon U_\infty \to X_\infty$. By going-down we can find $\beta \in U_\infty$ specializing to $\beta'$ and such that $f_\infty(\beta) = \alpha$. In particular, $u \coloneqq \beta(0)$ specializes to $\beta'$. Since $X$ is geometrically unibranch at $x$, we have $U$ is irreducible at $u$. Since $\beta$ specializes to $\beta'$ we have that $\beta(\eta)$ specializes to $\beta'(\eta)$. This implies that the only component of $U$ passing through $u$ is $U_i$ and since $U_i$ is smooth at $u'$, $U$ is smooth at $u$. Hence $X$ is smooth at $x = \alpha(0)$.
\end{proof}

In \cite[Question 2.10]{Reg18} Reguera asked whether, for each stable $\alpha \in X_\infty$ with $X$ being analytically irreducible at $\alpha(0)$, regularity of $\cO_{X_\infty,\alpha}$ implies that $X$ is nonsingular at $\alpha(0)$. As mentioned in the introduction, \cref{t:smooth-arcs} provides a stronger statement under the stronger assumption that $X$ is geometrically unibranch at $\alpha(0)$. Below we give an example where \cref{t:smooth-arcs} fails when assuming just unibranch instead. However, the local ring in question is not regular, meaning that Reguera's question is still open in the case of multiple geometric branches.

Let $X$ be the pinch point hypersurface singularity given $f = y^2 - x^2 z$ inside $\A^3_\C$. Take the blowup of $\A^3_\C$ with center the (reduced) singular locus of $X$ defined by the ideal $(x,y)$. In one affine chart, the morphism from the total transform $X'$ to $X$ is algebraically given by the ring homomorphism
\[
    \C[x,y,z]/(f) \to \C[x,u,v]/(x^2u), \: x \mapsto x, y \mapsto xv, z \mapsto v^2 - u.
\]
Note that the normalization $\widetilde{X}$ of $X$ is the closed subscheme of $X'$ cut out by $u = 0$. Now take the arc $\beta$ on $X'$ given by 
\[
    \C[x,u,v]/(u^2,v) \to \C(x_i,v_j \mid i\geq 1, j\geq 0)[[t]],\: x\mapsto \sum_{i\geq 1} x_i t^i, v \mapsto \sum_{j\geq0} v_jt^j, u \mapsto 0.
\]
Composing with $X' \to X$ gives an arc $\alpha$, which is the generic arc on $X$ having order of contact $1$ with $\Sing X$. In particular, we have that $x \coloneqq \alpha(0)$ is the generic point of $\Sing X$.

\begin{proposition}
    \label{p:unibranch-degen}
    Keep the notation from above. Then:
    \begin{enumerate}
        \item \label{eq:unibranch-degen1} $X$ is unibranch but not geometrically unibranch at $\alpha(0)$.
        \item \label{eq:unibranch-degen2} $\embcodim(\cO_{X_\infty,\alpha}) = 0$.
        \item \label{eq:unibranch-degen3} The local ring $\cO_{X_\infty,\alpha}$ is not reduced (and hence in particular not regular).
    \end{enumerate}
\end{proposition}

In fact, by e.g.\ \cite[Theorem 1.1]{BH21} the local ring $\cO_{X_\infty,\alpha}$ is not even Noetherian.

\begin{proof}
    For \eqref{eq:unibranch-degen1} we consider the preimage of $\alpha(0)$ under the normalization map, which by the above is easily seen to be a single point with quadratic residue field extension. To show \eqref{eq:unibranch-degen2}, we consider the projection $p \colon X \subset \A^3_\C \to \A^2_\C$ given by $(x,y,z) \mapsto (x,z)$. The image $\gamma= p_\infty(\alpha)$ is the generic arc with order of contact $1$ with the subscheme $\{x=0\}$. Therefore the local ring $\cO_{(\A^2_\C)_\infty,\gamma}$ is regular of dimension $1$. By \cite[Theorem 8.1]{CdFD22} we have $\embdim(\cO_{X_\infty,\alpha}) = 1$. Clearly $\dim(\cO_{X_\infty,\alpha}) \geq 1$, and by \cite[Theorem D]{CdFD24} the maximal ideal $\fm_\alpha$ of $\cO_{X_\infty,\alpha}$ is finitely generated. Now we may apply the argument in the proof of \cite[Corollary 4.6]{Reg06} to get $\dim(\widehat{\cO_{X_\infty,\alpha}}) = 1$.

    Finally, in order to show that $\cO_{X_\infty,\alpha}$ is not reduced we consider the following element of the coordinate ring of $(\A^3_\C)_\infty$:
    \[
        g \coloneqq 2y^{(0)}z^{(0)} x^{(1)} -2x^{(0)}z^{(0)}y^{(1)} + x^{(0)}y^{(0)}z^{(1)}.
     \]
    Then $g^2 = h_1 f^{(0)} + h_2 f^{(1)}$ with
    \[
h_1 \coloneqq 4(x^{(1)})^{2}(z^{(0)})^{2}
+2x^{(0)}x^{(1)}z^{(0)}z^{(1)}
-4(y^{(1)})^{2}z^{(0)}
+2y^{(0)}y^{(1)}z^{(1)}
\]
and
\[
h_2 \coloneqq -2x^{(0)}x^{(1)}(z^{(0)})^{2}
+2y^{(0)}y^{(1)}z^{(0)}
-(y^{(0)})^{2}z^{(1)}.
\]

    To show that the image of $g$ is nonzero in $\cO_{X_\infty,\alpha}$ it suffices to show that the image of $g$ in $\cO_{X'_\infty,\beta}$
    \[
        \widetilde{g} = 2(x^{(0)})^2u^{(0)}v^{(1)} - (x^{(0)})^2v^{(0)}u^{(1)}
    \]
    is nonzero. Write $k_\beta = k(x_i,v_j \mid i\geq1 , j\geq0)$ and $x(t) \coloneqq \sum_{i \geq 1} x_i t^i$, $v(t) \coloneqq \sum_{j\geq0} v_j t^j$. Let $R$ be the valuation ring of $k_\beta(r,s)$ with respect to the rank $2$-valuation $v(r) = \varepsilon$, $v(s) = 1$ with value group $\Gamma = \Z \oplus \varepsilon \Z$ endowed with the lexicographic order, as in \cite[Example 5.12]{CdFD24}. Consider the infinitesimal deformation $k[x,u,v]/(x^2u) \to R/(r^2s)[[t]]$ of $\beta$ given by
    \[
        x \mapsto r - x(t), v \mapsto v(t), u \mapsto s \left(1 + \frac{1}{r} x(t) + \frac{1}{r^2} x(t)^2 + \ldots \right)^2.
    \]
    The image of $\widetilde{g}$ under this map equals $2x_1v_0rs \neq 0$ and hence $\widetilde{g} \neq 0$ in $\cO_{X'_\infty,\beta}$.
\end{proof}

\begin{remark}
    The particular nilpotent $g$ was found as follows: the surface $X$ is isomorphic to the deformation of the nodal to the cuspidal plane curve singularity. Thus we can use the strategy for finding nilpotents in the arc space of plane curve singularities in \cite{Seb11} and then deform. The suspicion that this nilpotent does not vanish in the local ring at $\alpha$ came from \cite[Example 5.13]{CdFD24}.
\end{remark}

\begin{remark}
    \label{r:reduced-main-cpt}
    Note that the pair $(\A^3,X)$ is log canonical, but the singularities of $X$ are not rational. By \cite[Example 7.12]{EM09} all jet schemes $X_m$ are equidimensional but not irreducible. Let $X_m^0$ be the unique irreducible component of $X_m$ which contains $\pi_{m,0}^{-1}(X \setminus \Sing X)$, where $\pi_{m,0} \colon X_m \to X$. We consider $X_m^0$ endowed with its natural (a priori nonreduced) scheme structure and call it the \emph{main component} of $X_m$. Since $X\setminus \Sing X$ is smooth, we have that $X_m^0$ is irreducible and generically reduced. Since $X$ is a hypersurface it follows from the unmixedness theorem that $X_m^0$ is reduced. Note that the nilpotent element $g$ indeed vanishes on the main component of $X_1$, while it does not vanish in the local ring at $\alpha$, whose image in $X_1$ is not contained in $X_1^0$.
\end{remark}

\printbibliography

\end{document}